\documentclass[11pt,a4paper]{article}

\usepackage{amsmath,amsthm,amssymb}
\usepackage{times}
\usepackage{enumerate}
\usepackage[textwidth=160mm, top=30mm, bottom=30mm,hmarginratio=1:1]{geometry}

\usepackage{setspace}

\usepackage[mathlines]{lineno}

\usepackage[dvipsnames]{xcolor}

\newcommand*\patchAmsMathEnvironmentForLineno[1]{%
  \expandafter\let\csname old#1\expandafter\endcsname\csname #1\endcsname
  \expandafter\let\csname oldend#1\expandafter\endcsname\csname end#1\endcsname
  \renewenvironment{#1}%
     {\linenomath\csname old#1\endcsname}%
     {\csname oldend#1\endcsname\endlinenomath}}%

\newcommand*\patchBothAmsMathEnvironmentsForLineno[1]{%
  \patchAmsMathEnvironmentForLineno{#1}%
  \patchAmsMathEnvironmentForLineno{#1*}}%

\AtBeginDocument{%
  \patchBothAmsMathEnvironmentsForLineno{equation}%
  \patchBothAmsMathEnvironmentsForLineno{align}%
  \patchBothAmsMathEnvironmentsForLineno{flalign}%
  \patchBothAmsMathEnvironmentsForLineno{alignat}%
  \patchBothAmsMathEnvironmentsForLineno{gather}%
  \patchBothAmsMathEnvironmentsForLineno{multline}%
}

\usepackage{dsfont}
\DeclareMathOperator{\Int}{Int}

\DeclareMathOperator{\tr}{Tr}

\DeclareMathOperator{\interior}{int}
\DeclareMathOperator{\closure}{cl}

\newcommand{\R}{\mathbb R}

\newcommand{\gep}{\varepsilon}
\newcommand{\RR}{\mathbb{R}}
\newcommand{\RRn}{\mathbb{R}^n}

\newcommand{\open}{(\Omega_t)_{t\in(0,T)}}

\newcommand{\omegaone}{(\Omega^1_t)_{t\in(0,T)}}
\newcommand{\omegatwo}{(\Omega^2_t)_{t\in(0,T)}}
\newcommand{\close}{(\mathcal{F}_t)_{t\in(0,T)}}

\newcommand{\limsupstard}{{\limsup_{\varepsilon\to0^+}}^*\;}
\newcommand{\liminfstard}{\liminf_{\varepsilon\to0^+}\hspace{-0.08cm}{}_*\;}

\newcommand{\omegaeb}{\omega^{\varepsilon,\beta}}
\providecommand{\norm}[1]{\lVert#1\rVert}
\providecommand{\abs}[1]{\lvert#1\rvert}

\theoremstyle{definition}
\newtheorem{thm}{Theorem}[section]
\newtheorem{cor}[thm]{Corollary}

\newtheorem{lem}[thm]{Lemma}
\newtheorem{defin}[thm]{Definition}
\newtheorem{rem}[thm]{Remark}

\newcommand{\subjclass}[1]{\bigskip\noindent\emph{2010 Mathematics Subject Classification:}\enspace#1}
\newcommand{\keywords}[1]{\noindent\emph{Keywords:}\enspace#1}

\numberwithin{equation}{section}

\begin{document}

\title{Singular limits of the subelliptic Allen-Cahn equation\\ and geometric flows}

\author{Cecilia De Zan\thanks{current address: Liceo Statale "Giuseppe Berto", Ministero Istruzione, Universit\`a e ricerca, Mogliano Veneto (Italy). email: dezancecilia@gmail.com} \;\& Pierpaolo Soravia\thanks{email: soravia@math.unipd.it. Supported in part by Universit\`a di Padova research project 2025DM1DOR-00669 "Equazioni differenziali nonlineari, problemi variazionali e di controllo ottimale" }\\
Dipartimento di Matematica\\ Universit\`{a} di Padova, via Trieste 63, 35121 Padova, Italy}

\date{}
\maketitle

\begin{abstract}
This paper investigates the singular limits of a bistable reaction-diffusion equation in a subelliptic setting, where the standard Euclidean diffusion is replaced by a subelliptic operator defined by a family of vector fields. We introduce a Riemannian approximation of these vector fields to regularize the problem, yielding a uniformly parabolic equation whose ellipticity degenerates as the singular parameter vanishes.

We prove that the family of solutions converges globally in time to the stable equilibria away from a propagating interface. Our approach uses the definition of generalized geometric flow developed by Barles--Souganidis, uniquely defined through the level-set equation and viscosity solutions. The behavior of the interface depends on the traveling wave speed $c$: when $c \neq 0$, the front propagates with constant normal velocity governed by a first-order subelliptic eikonal equation.

When $c = 0$, the front evolves according to the horizontal mean curvature flow equation defined by the original vector fields. At present, this case comes with limitations: it requires a more restrictive condition on test functions at characteristic points to define generalized flows, where the horizontal gradient vanishes and standard arguments break down. This framework relies on our earlier work~\cite{deso5}, where the equivalence between the generalized flow and the level-set approach was established when the vector fields generate a step-two Carnot group. We also assume that the comparison principle holds for the level-set equation in order for the level-set approach to be well-posed.
\end{abstract}


\section{Introduction}

When $f:\R\to\R$ is a cubic nonlinearity with exactly three zeroes $m_-<m_o<m_+$, with
$m_\pm$ being the stable equilibria and $m_o$ the unstable one, reaction-diffusion equations of the form
\begin{equation}\label{reac-diff eq no scaling}
\phi_t(x,t)-\Delta\phi(x,t)+f(\phi(x,t))=0,\quad(x,t)\in\RRn\times(0,+\infty)
\end{equation}
arise naturally in a variety of mathematical models, such as phase transitions and flame propagation. This equation was originally introduced by Allen and Cahn \cite{ac} to model the motion of antiphase boundaries in crystalline solids.
In most of these applications, one can observe for large times the development of fronts as the boundaries of the regions where the solution $\phi$ converges to the stable equilibria of the cubic $f$.
The connection between the Allen-Cahn equation and the geometry of interfaces is also reflected at the variational level, as initially shown by the $\Gamma-$convergence result of Modica, Mortola \cite{momo}. Caffarelli, Cordoba \cite{cc} established uniform convergence of the singular perturbation problem, providing quantitative estimates on the thickness of the transition layer separating the two stable phases. Recently Chan, Fernández-Real, Figalli, Serra \cite{cffs}
 proved that global stable solutions to the free-boundary Allen--Cahn equation in three dimensions are one-dimensional, confirming a variant of the De Giorgi conjecture.
 It is worth mentioning that the analogous one dimensional symmetry property is known to fail in the sub-Riemannian setting. Birindelli, Lanconelli \cite{bl} exhibited a counterexample of stable but not one dimensional solution to the De Giorgi conjecture in the Heisenberg group. This result has been further investigated by Ferrari, Valdinoci \cite{fv} and Birindelli, Ferrari, Valdinoci \cite{bfv}. This illustrates how the subelliptic setting may give rise to qualitatively different phenomena with respect to the Euclidean case.

From the parabolic point of view, when
\begin{equation}\label{f example}
f(q)=2(q-m_-)(q-m_o)(q-m_+),
\end{equation}
formal results of Fife \cite{fi,fi2} and Caginalp 
\cite{cagi1, cagi2, cagi3} show that the fronts propagate with normal velocity
\begin{equation}\label{vel}
V=c+\frac{\kappa}{t}+O\Big(\frac{1}{t^2}\Big),\qquad t\gg1,
\end{equation}
where $c=2m_o-m_+-m_-$, and $\kappa$ denotes the mean curvature of the front. 
In more general cases $c$ comes from solving the traveling wave equation.
A rigorous justification of this result was later given by 
Evans, Soner, Souganidis \cite{ess} and Barles, Soner, Souganidis \cite{bss}, 
in the more general setting where the cubic function depends also on 
the variables $x,t$ (i.e.\ $f=f(q,x,t)$) and the second order operator 
is anisotropic but uniformly elliptic.
To obtain the first term in (\ref{vel}) one 
considers the singular limit as $\gep\to0^+$ of the rescaled equation in 
(\ref{reac-diff eq no scaling}) with the scaling $(\gep^{-1}x,\gep^{-1}t)$. 
If $c=0$ one 
must proceed to the next time-scaling $(\gep^{-1}x,\gep^{-2}t)$ to obtain the second term. This requires analyzing two singular equations:
$$\phi^\gep_t(x,t)-\gep\Delta\phi^\gep(x,t)+\frac{f(\phi^\gep(x,t))}{\gep}=0,$$
when $c\neq0$, and
$$\phi^\gep_t(x,t)-\Delta\phi^\gep(x,t)+\frac{f(\phi^\gep(x,t))}{\gep^2}=0,$$
when $c=0$.

In this paper, we extend this analysis to the corresponding subelliptic setting, replacing standard Euclidean derivatives with derivatives along a family of $m$ vector fields $X_1, \dots, X_m$, with $m < n$. Although we are not aware of direct applications employing the subelliptic Allen-Cahn equation itself, the resulting limiting equation, the sub-Riemannian mean curvature flow equation, is of significant interest in neurobiology and visual cortex modeling (see, e.g., Citti, Manfredini, and Sarti~\cite{cms}, Baspinar and Citti~\cite{baci}, and the references therein).
In these models, the underlying geometry is non-Euclidean. Signals propagate preferentially along specific directions, naturally described by a family of vector fields rather than the full gradient. The resulting sub-Riemannian structure captures this anisotropy in a mathematically precise way, while the degeneracy of the diffusion matrix reflects the constraint that diffusion occurs exclusively along the directions spanned by the vector fields. This inherent degeneracy renders the analysis significantly more delicate than in the classical setting.
We consider the equation
\begin{equation}\label{eqnorescaling}
u_t(x,t)-\tr(X^2u(x,t))+f(u(x,t))=0,\quad\mbox{ in }\RRn\times(0,+\infty).
\end{equation}
The notation $Xu$ indicates the horizontal gradient with respect to the 
vector fields $\{X_1,\dots,X_m\}$, seen as differential operators,
\begin{equation}\label{eqfields}
X_j=\sum_{i=1}^n\sigma_{ij}(x)\partial_i,\quad j\in\{1,\dots,m\},
\end{equation}
i.e. $Xu(x)=\nabla u(x)\sigma(x)$ if $\sigma$ is the matrix whose columns are the vectors $\sigma^j=(\sigma_{1j},\dots,\sigma_{nj})^T$,
while $X^2u$ indicates the symmetrization of the corresponding horizontal Hessian $(X_k(X_j u))_{j,k\in\{1,\dots,m\}}$.

To overcome the fact that the equation is not uniformly parabolic, we introduce a Riemannian approximation
 obtained by augmenting the original vector fields with small isotropic corrections of order $\gep^k$. This regularization renders the approximating equations uniformly parabolic for each $\gep>0$, while recovering the original degenerate subelliptic structure in the singular limit. The choice of the exponent $k$ depends on the two cases below and plays a crucial role in the estimates, particularly near characteristic points, to control remainder terms. It will be carefully calibrated in the analysis, although we do not seek an optimal choice.
We therefore analize the singular limits as $\gep \to 0^+$ of two approximating problems. 

\medskip
\noindent\textbf{Case 1} ($c \neq 0$):
\begin{equation}\label{case1}
u^\gep_t(x, t) 
- \gep\,\tr(Y^2_\gep u^\gep(x, t)) 
+ \frac{f(u^\gep(x, t))}{\gep} = 0 
\quad \mbox{in } \RRn \times (0, +\infty).
\end{equation}

\noindent\textbf{Case 2} ($c = 0$):
\begin{equation}\label{case2}
u^\gep_t(x, t) 
- \tr(Y^2_\gep u^\gep(x, t)) 
+ \frac{f(u^\gep(x, t))}{\gep^2} = 0 
\quad \mbox{in } \RRn \times (0, +\infty).
\end{equation}

\medskip
\noindent 
Here, $Y^2_\gep u^\gep$ denotes the symmetrized horizontal Hessian associated with a Riemannian approximation of the family $\{X_i\}$. Notice that the approximation parameter is the singular parameter of the problem.
When $c\neq0$, the dominant effect in the limit is the transport of the front at constant speed, and the interface dynamics are governed by the first-order Hamilton-Jacobi eikonal equation
    \begin{equation}\label{asymp eq}
        w_t(x, t) + c\,|X w(x, t)| = 0, 
        \quad (x, t) \in \RRn \times (0, +\infty),
    \end{equation}
    where the front is identified as a level set of $w$. Note that the Hamiltonian is not coercive as it only depends on the horizontal gradient.
    Our result in this case holds true for a general family of vector fields.
In contrast, when $c=0$, the first-order contribution vanishes and one must proceed to a finer time-scale to capture the front motion, which turns out to be driven by the geometry of the interface itself through the horizontal mean curvature. 
Our analysis is here limited to step two Carnot groups.
The front evolves according to the equation
 \begin{equation}\label{allen cahn eps eq}
 \begin{aligned}
 w_t(x,t)-\sum_{i,j=1}^m\left(\delta_{ij}-
 \frac{X_iw(x,t)X_jw(x,t)}{\sum_{k=1}^m
 (X_kw(x,t))^2}\right)X_iX_jw(x,t)=0,
 \end{aligned}
 \quad (x, t) \in \RRn \times (0, +\infty),
 \end{equation}
which is the horizontal mean curvature flow equation  defined by the original $m$ vector fields.
We will assume that equation \eqref{allen cahn eps eq} satisfies the comparison principle for viscosity solutions in order to carry out the analysis in this case with a well posed level-set approach. This fact is not completely solved in the current literature as we discuss later in this section.

Starting from initial data assuming values on both sides of $m_o$, we prove that the family of solutions 
$u^\gep$ converges globally in time to the stable equilibria $m_\pm$ away from a propagating interface (front). 
In both cases, global in time convergence is obtained via the abstract geometric framework of Barles and Souganidis \cite{bs}, using generalized flows, which is designed to study singular limits leading to moving interfaces.  When $c\neq0$, the equivalence between the level-set approach and generalized flows relies on standard comparison results, which hold unchanged in our setting.
    In the case $c=0$, singularities occur not only where the gradient of the solution vanishes, but at all characteristic points, i.e. when the horizontal gradient vanishes.
This more complex singular structure, breaks standard arguments. To address this, besides using a Riemannian approximation, we enforce a stricter condition on the test functions at characteristic points within our definition of generalized flow, recovering regularity to carry out the asymptotic analysis.
Consequently, our framework is currently restricted to step-two Carnot groups, where the equivalence between generalized flows and the level-set approach was established in \cite{deso5}. We note, however, that any future extension of that equivalence would immediately broaden our result.
For $c=0$, the proof of the mentioned equivalence also requires that we use in the subelliptic mean curvature flow equation an alternative  definition of viscosity solution due to Barles, Georgelin \cite{bage}, whose extension to the sub-Riemannian framework is again presently limited to step two Carnot groups, see e.g. Ferrari, Liu, Manfredi \cite{flm} and the authors \cite{deso5}.
Finally, since as we mentioned, the comparison principle for the level-set equation remains open in full generality, global convergence for $c=0$ highlights key open problems in the general theory of subelliptic mean curvature flows.

We recall that the geometric approach in \cite{bs} bypasses the need for smoothness assumptions on both the interface and the solution of the geometric equation, common in the formal arguments, allowing the asymptotic analysis to be carried out in a more direct and robust way, see also Souganidis \cite{soug}. Our analysis follows the approach of \cite{bs}, as revisited by Barles and Da Lio \cite{bdl} in the context of bounded domains with Neumann boundary conditions, and we frequently combine arguments from both works. We show that it extends naturally to our subelliptic framework.
We have to tackle the possible presence of characteristic points, where the normal to the front is well defined but can be orthogonal to all the vector fields of the family.
We thus face new difficulties in the construction of global super- or sub-solutions of the reaction diffusion equation (the {\it initialization step}) due to degeneracy of the diffusion term. We handle them by linking the degeneracy to the singular parameter $\gep\to0$. We face similar difficulties also in the proof that the asymptotics defines a generalized flow (the {\it propagation step}).
To handle general initial conditions, we also draw on ideas from Chen \cite{xce} to establish that an interface forms within a very short time.
In \cite{dzs2}, we utilized similar ideas to address the singular limits when the front exhibits a discontinuous velocity, a rather different framework.
Finally, we mention the work of Da Lio, Kim, Slepčev \cite{dlks}, on the asymptotics of a reaction-diffusion equation with a nonlocal term in a bounded domain with a nonlinear oblique derivative boundary condition.

We briefly review the relevant literature on the level set approach, geometric equations and subelliptic mean curvature flows. The level set method was proposed by Osher, Sethian \cite{os} for numerical computations of moving fronts.
The rigorous theory of weak front evolution started with the work by Evans, Spruck \cite{es} for the mean curvature flow and by Chen, Giga, Goto \cite{cgg} for more general geometric equations. The reader is also referred to the book by Giga \cite{gi}, where the approach is discussed in detail, 
see also Souganidis \cite{soug} and the references therein 
for the main applications of the theory.
The study of horizontal normal and horizontal mean curvature was pioneered by Danielli, Garofalo, and Nhieu \cite{dagani}. Since then, the resulting subelliptic geometric flows have drawn considerable attention. Capogna and Citti \cite{caci} established existence in Carnot groups via vanishing viscosity Riemannian approximations, while Dirr, Dragoni, and Von Renesse \cite{didr} utilized stochastic approximations for more general Hörmander structures. 
The comparison principle only for axisymmetric fronts is also proved in Ferrari, Liu, Manfredi \cite{flm}.
Uniform regularity estimates for the flow of graphs have also been proven by Capogna, Citti, and Manfredini \cite{cacima}. However, uniqueness and comparison principles remain challenging. While \cite{caci} proved the comparison principle under specific initial conditions to avoid characteristic points, recently Baspinar and Citti \cite{baci} established the comparison principle in step-two Carnot groups and in the 3d rototranslation group. They prove comparison principle between solutions that are limits of suitable Riemannian regularizations. We note that \cite{caci} and \cite{baci} employ notions of solutions that differ from standard viscosity solutions at characteristic points but are compatible with ours.

The paper is organized as follows. Section 2 introduces the main definitions and notation, including the subelliptic framework, a Riemannian approximation of the diffusion matrix, and the notion of generalized super- and subflows adapted to the sub-Riemannian setting. In particular, we extend the classical definition of Barles, Souganidis \cite{bs} by adding an extra condition at characteristic points, and we describe the abstract method used to study the singular limits. Section 3 is devoted to the first asymptotic problem $c\neq0$. We prove that the family of solutions converges to the stable equilibria away from a front propagating with constant normal velocity c.
Section 4 addresses the second asymptotic problem $c=0$. Under the additional assumption that the vector fields generate a step-two Carnot group, and that the flow equation satisfies the comparison principle, we prove convergence to a front evolving by horizontal mean curvature. Finally, the Appendix collects auxiliary known results on the traveling wave equation, the structure of step-two Carnot groups, and the equivalence between Barles-Souganidis generalized flows and the level-set approach.

\section{Preliminaries on degenerate Allen Cahn equations and generalized flows}

We introduce first our notation.
Let $X_j=\sum_{i=1}^n\sigma_{ij}(x)\partial_i$, $j=1,\dots,m$, $m<n$, be a 
family of vector fields with coefficients the column vectors 
$\sigma^{j}=(\sigma_{ij})_{i=1,\dots,n}$,
$\sigma^{j}\in C^2(\RRn,\RRn)\cap W^{2,\infty}(\RRn,\RRn)$.
We define the so-called \emph{horizontal gradient} of a $C^1$-function $\Phi$ as
$$X\Phi(x):=(\langle \sigma^{j}, \nabla \Phi(x)\rangle)_{j=1,\dots,m}=\nabla\Phi(x)\sigma(x),\qquad x\in\RRn,$$
where $\sigma=(\sigma_{ik})_{i,k}\in C^2(\RRn,\RR^{n\times m})\cap 
W^{2,\infty}(\RRn,\RR^{n\times m})$ is the matrix map whose columns are $\sigma^{1},\dots,\sigma^{m}$.
We also define the \emph{symmetrized horizontal Hessian matrix} of $\Phi$, $X^2\Phi(x)$, 
by putting first
\begin{equation}\label{eqdrift}
\begin{aligned}
\tilde X^2\Phi(x) &:= (X_jX_i\Phi(x))_{i,j=1,\dots,m} = \Big(\sigma^{j}\cdot 
                     \nabla(\sigma^{i}\cdot\nabla\Phi(x))\Big)_{i,j=1,\dots,m}\\
                  &= \sigma^T(x)D^2\Phi(x)\sigma(x)
                     +\left(\nabla\Phi(x)D(\sigma(x))\right)\sigma(x),
\end{aligned}
\end{equation}
where the first order non-symmetric drift term appears
\begin{equation*}
\left(\nabla\Phi(x)D(\sigma(x))\right)\sigma(x)=
\Big((\nabla\Phi(x)D\sigma^{i}(x))\sigma^{j}(x)\Big)_{i,j=1,\dots,m}=\sum_{l,k=1,\dots,n}\partial_lu\;\partial_k(\sigma_{li})\;\sigma_{kj},
\end{equation*}
explicitly written in coordinates. Next
$$X^2\Phi(x)=(\tilde X^2\Phi(x))^*=
\frac{\tilde X^2\Phi(x)+\tilde X^2\Phi(x)^T}{2}.$$
We also define a matrix map 
$A(\cdot)\in C^2(\RRn,\mathcal{S}^n)\cap W^{2,\infty}(\RRn,\mathcal{S}^n)$ as
$$A(x)=\sigma(x)\sigma^T(x),\quad\mbox{ for any }x\in\RRn.$$
Clearly $A$ satisfies the following conditions,
\begin{eqnarray}\label{eqcore1}
&&\mbox{for all }i,j,k\in\{1,\dots,n\},\,a_{ij},a_{ij,x_k}
\mbox{ are bounded and continuous on }\RRn,\label{A regularity}\\
&&A(x)q\cdot q\geq0\mbox{ for any }(x,q)\in\RRn\times\RRn.\label{A non negative}
\end{eqnarray}
Moreover taking the trace
$$\tr(X^2\Phi(x))=\tr(\tilde X^2\Phi(x))=\tr(A(x)D^2\Phi(x))
+\tr(\nabla\Phi(x)D(\sigma(x))\sigma(x)).$$

We study the asymptotic behavior of the solutions $u$ of the 
Cauchy problems for appropriate rescalings of the equation \eqref{eqnorescaling}
when the map $f\in C^2(\R)$ is a cubic-type nonlinearity such that
\begin{equation}\label{f property}
\begin{array}{l}
f\mbox{ has exactly three zeroes }m_-<m_o<m_+,\\
f(s)>0\mbox{ in }(m_-,m_o)\cup(m_+,+\infty)\mbox{ and }
f(s)<0\mbox{ in }(-\infty,m_-)\cup(m_o,m_+),\\
f'(m_\pm)\geq d>0,\,f'(m_o)<-d<0,\,f''(m_-)<0\mbox{ and }f''(m_+)>0.
\end{array}
\end{equation}

Writing in coordinates $X^2u$, the equation also reads
\begin{equation}\label{asymp eq deg}
u_t(x,t)-\tr(A(x) D^2u(x,t))-\tr((\nabla u(x)D\sigma(x))\sigma(x))+f(u)=0.
\end{equation}
The matrix $A(x)$ is only positive semi-definite, as $m<n$, and equation (\ref{asymp eq deg}) 
does not satisfy a uniform ellipticity condition. 
We introduce a Riemannian approximation of the matrix $A$ into the scaling.
We start by enlarging the family of the vector fields and define a new matrix map 
$\sigma_\gep\in C^2(\RRn,\RR^{n\times (m+n)})\cap 
W^{2,\infty}(\RRn,\RR^{n\times (m+n)})$, $\gep>0$ as
\begin{equation}\label{sigmaeps def}
\sigma_\gep(\cdot)=\left(\sigma(\cdot)\;|\;\gep^k\, I_n\right),
\end{equation}
where for now $0<k<1$ and $I_n$ denotes the $n\times n$ identity matrix. 
As Riemannian approximation of $A$ we then have
$$A_\gep(\cdot)=\sigma_\gep(\cdot)\sigma_\gep^T(\cdot)=A(\cdot)+\gep^{2k} I_n.$$
Moreover $A_\gep\to A$, uniformly as $\gep\to0^+$.
For any $x\in\RRn$ the matrix $A_\gep(x)\in\RR^{n\times n}$ is therefore uniformly positive 
definite and the ellipticity constant degenerates as $\gep\to0^+$.
Indeed
\begin{equation}\label{Agep pos def}
A_\gep(x)q\cdot q= A(x)q\cdot q+\gep^{2k}\abs{q}^2\geq\gep^{2k}\abs{q}^2.
\end{equation}
If we introduce the new family of vector fields by putting $Y_{\gep j}=X_j$ for $j=1,\dots,m$ and $Y_{\gep j}=\gep^k\partial_j$ if $j=m+1,\dots,m+n$ then, for smooth $\Phi$,
$$Y_\gep \Phi(x)=(\langle\nabla \Phi(x), \sigma_\gep^j(x)\rangle)_{j=1,\dots,m+n}=\nabla \Phi\sigma_\gep=
\left(X\Phi(x)\;|\;\gep^k\nabla \Phi(x)\right),$$
$$\tilde Y_{\gep}^2\Phi(x)=
\left(\begin{array}{cc}
\tilde X^2\Phi(x)&\gep^kD(\nabla \Phi(x)\sigma(x))\\
\gep^kD^2\Phi(x)\sigma(x)&\gep^{2k}D^2\Phi(x)
\end{array}\right),\qquad Y_{\gep}^2\Phi(x)=(\tilde Y_\gep^2\Phi(x))^*.$$
Concerning the drift term in \eqref{eqdrift}, observe that
$$\left(\nabla \Phi(x)D\sigma_\gep(x)\right)\sigma_\gep(x)
=\left(\begin{array}{cc}
\left(\nabla \Phi(x)D\sigma(x)\right)\sigma(x)&\gep^k\nabla \Phi(x)D\sigma(x)\\
0_{n\times m}&0_{n\times n}
\end{array}\right),
$$
since the last $n$ columns of $\sigma_\gep$ do not depend on the variables. In particular
\begin{equation}\label{drift disappears}
\hbox{Tr}((\nabla \Phi(x)D\sigma_\gep(x))\sigma_\gep(x))=\hbox{Tr}((\nabla \Phi(x)D\sigma(x))\sigma(x))
\end{equation}
and we are doing a vanishing viscosity approximation to the parabolic equation as 
\begin{equation}\label{approximating trace}
\hbox{Tr}(Y_\gep \Phi)=\hbox{Tr}(X \Phi)+\gep^{2k}\hbox{Tr}(D^2 \Phi).
\end{equation}
\begin{rem}\label{remgradest} We observe for later convenience that if
$\Phi\in C^1(\RRn)$, $x\in\RRn$ and $\gep>0$,
$$\abs{Y_\gep \Phi(x)}^2=\abs{X\Phi(x)}^2+\gep^{2k}\abs{\nabla \Phi(x)}^2\geq\gep^{2k}\abs{\nabla \Phi(x)}^2,$$
and 
$\abs{Y_\gep \Phi(x)}\longrightarrow\abs{X \Phi(x)}$, locally uniformly as $\gep\to0^+$. Furthermore
$$\{x\in\RRn:X\Phi(x)=0\}\supseteq\{x\in\RRn:\nabla \Phi(x)=0\}=\{x\in\RRn:Y_\gep \Phi(x)=0\}.$$
In the study of geometric front evolutions $t\to\{\phi(x,t)=0\}$, one frequently encounters points where the Euclidean normal is well-defined ($\nabla\phi(x,t)\neq0$), yet the point is \emph{characteristic} with respect to the underlying geometry, meaning $X\phi(x,t)=\nabla\phi(x,t)\sigma(x)=0$. In the subset of characteristic points, the regularized horizontal gradient satisfies
$$0\neq |Y_\gep\phi(x,t)|=\gep^k\abs{\nabla\phi(x,t)}\xrightarrow[\gep\to0^+]{}0,$$
locally uniformly.
Consequently, when analyzing the asymptotic behavior of equations \eqref{asymp eq} and \eqref{allen cahn eps eq}, the choice of the scaling exponent k in a Riemannian approximation is crucial. It must be precisely calibrated to balance the vanishing regularizing effect with the underlying subelliptic degeneracy near the characteristic set.
\end{rem}

\subsection{Generalized super and subflows}

Equations \eqref{asymp eq} and \eqref{allen cahn eps eq} are examples of geometric equations describing flows of interfaces. In general, let
$$
F : \mathbb{R}^n 
\times \mathbb{R}^m \setminus \{0\} 
\times \mathcal{S}^m \to \mathbb{R}
$$
be a locally bounded and continuous function, where $\mathcal{S}^m$ denotes 
the space of $m \times m$ symmetric matrices. We impose the following 
structural conditions on $F$.
\begin{description}
  \item[(F1)] \textit{Consistency at the singularity:}
    \begin{equation*}
      F^*(x,0,\mathbb{O}) = F_*(x,0,\mathbb{O}), 
      \quad \text{for all } x \in \mathbb{R}^n ;
    \end{equation*}
  \item[(F2)] \textit{Ellipticity:} for any 
    $x \in \mathbb{R}^n$ 
    , $p \in \mathbb{R}^m \setminus \{0\}$, and $A, B \in \mathcal{S}^m$,
    \begin{equation*}
      F(x,p,A) \leq F(x,p,B) \quad \text{if } A \geq B;
    \end{equation*}
  \item[(F3)] \textit{Geometric condition:} for all $\lambda > 0$, 
    $\mu \in \mathbb{R}$, $x \in \mathbb{R}^n$ 
    , $p \in \mathbb{R}^m \setminus \{0\}$, and $A \in \mathcal{S}^m$,
    \begin{equation*}
      F(x,\lambda p, \lambda A + \mu(p \otimes p)) 
      = \lambda\, F(x,p,A).
    \end{equation*}
\end{description}
Here, the lower and upper semicontinuous extensions of $F$ at singular 
points are defined by
$$
F_*(x,0,A) = \lim_{r \to 0^+} \inf\{F(y,q,B) : q \neq 0,\; 
|(y,q,B)-(x,0,A)| \leq r\}
$$
and similarly to define $F^*$ by replacing $\inf$ with $\sup$. In particular, 
conditions (F1)--(F3) together imply
$$
F_*(x,0,\mathbb{O}) = 0, \quad \text{for all } 
x \in \mathbb{R}^n.$$
Consider the {\it geometric} subelliptic equation 
\begin{equation}\label{eqlevelset}
v_t(x,t) + F(x,\, Xv,\, X^2 v) = 0, 
\quad x \in \mathbb{R}^n 
\end{equation}
in which only the horizontal first and second derivatives of $u$ appear. 
We are here only interested in two specific cases of the equation (\ref{eqlevelset}), either when $F:\R^m\to\R$ and $F(p)=c|p|$, $c\neq0$, or $F(p,A)=\hbox{Tr}\left((I_m+\frac{p}{|p|}\otimes \frac{p}{|p|})A\right)$. 
Therefore, in standard Euclidean coordinates of $\mathbb{R}^n$, the operator in~(\ref{eqlevelset}) 
takes the form
\begin{equation}\label{eqopg}
G(x,p,A) = F\!\left(x,\, p\,\sigma(x),\; 
\sigma^T(x)\, A\, \sigma(x)+(pD\sigma(x))\sigma(x)\right),
\end{equation}
where
$$
G : \Bigl(\mathbb{R}^n 
\times \mathbb{R}^n \setminus \{(x,p) : p\,\sigma(x) = 0\}\Bigr) 
\times \mathcal{S}^n \to \mathbb{R}.
$$
\begin{rem}
We show now that the operator $G$ inherits the geometric properties of $F$.
However, 
its singularity set is no longer just the origin, but the entire subset
$$
S = \{(x,p,A) \in \mathbb{R}^n 
\times \mathbb{R}^n 
\times \mathcal{S}^n : p\,\sigma(x) = 0\},
$$
where, for fixed $(x,A)$, the set $\{p : (x,p,A) \in S\}$ is a varying 
subspace that is not trivial since $m<n$ and $\mathcal{X}$ does not span 
$\mathbb{R}^n$ at $x$. In this sense, in the second instance above, $G$ falls outside the framework of the classical theory of the level set method. Indeed, ellipticity of $G$ follows 
since if $A \geq B$, then $\sigma^T(x) A \sigma(x)
\geq \sigma^T(x) B \sigma(x)$, and thus $G(x,t,p,A) \leq G(x,t,p,B)$. 
The geometric condition for $G$ follows from
$$
\begin{aligned}
G(x,\lambda p,\, \lambda A + \mu(p \otimes p))
&= F\!\left(x,\, \lambda p\,\sigma(x),\; 
\lambda\, \sigma^T(x) A \sigma(x) +\lambda(pD\sigma(x))\sigma(x)
+ \mu\,(p\,\sigma(x) \otimes p\,\sigma(x))\right) \\
&= \lambda\, F\!\left(x,\, p\,\sigma(x),\; 
\sigma^T(x) A \sigma(x)+(pD\sigma(x))\sigma(x)\right) 
= \lambda\, G(x,p,A),
\end{aligned}
$$
so (F2) and (F3) hold for $G$ as well.
Moreover, the lower semicontinuous extension of $G$ at the singular points 
where $p\,\sigma(x) = 0$ is given by
$$
\begin{aligned}
G_*(x,p,A) 
&= \lim_{r \to 0^+} \inf\bigl\{G(y,q,B) : 
(y,q,B) \in \operatorname{dom}(G),\; 
|(x-y, p-q, A-B)| \leq r\bigr\} \\
&= \lim_{r \to 0^+} \inf\bigl\{F(\, q\,\sigma(y),\, 
\sigma^T(y) B \sigma(y)+(q\;D\sigma(y))\sigma(y)),\\
&\qquad q\,\sigma(y) \neq 0,\; 
|(x-y, t-s, p-q, A-B)| \leq r\bigr\}
\end{aligned}
$$
and the following chain of inequalities holds:
$$\begin{array}{ll}
F_*\!\left(x,0,\, \sigma^T(x) A \sigma(x)+(pD\sigma(x))\sigma(x)\right) 
&\leq G_*(x,p,A) \leq G^*(x,p,A) \\
&\leq F^*\!\left(x,0,\, \sigma^T(x) A \sigma(x)+(pD\sigma(x))\sigma(x)\right).
\end{array}$$
Consequently, if also $\sigma^T(x) A \sigma(x) +(pD\sigma(x))\sigma(x)= \mathbb{O}$, then 
$G_*(x,p,A) = G^*(x,p,A) = 0$, and the following analogue 
of (F1) holds for $G$:
\begin{equation}\label{G0 condition}
G^*(x,p,A) = G_*(x,p,A) = 0,
\end{equation}
for all $(x,p,A) \in \mathbb{R}^n \times \mathbb{R}^n \times 
\mathcal{S}^m$ such that $p\,\sigma(x) = 0$ and $\sigma^T(x) A \sigma(x)+(pD\sigma(x))\sigma(x) = \mathbb{O}$. This is a stronger version of (F1).
\end{rem}

We recall the notion of viscosity solution for equation~(\ref{eqlevelset}).
\begin{defin}\label{defviscsol}
An upper (respectively lower) semicontinuous function 
$v : \mathbb{R}^n \times (0,+\infty) \to \mathbb{R}$ is a 
\emph{viscosity subsolution} (respectively \emph{supersolution}) 
of~(\ref{eqlevelset}) if and only if for any 
$\phi \in C^2(\mathbb{R}^n \times (0,+\infty))$, whenever 
$(x,t) \in \mathbb{R}^n \times (0,+\infty)$ is a local maximum 
(respectively minimum) point of $v - \phi$, we have
\begin{equation}\label{eqsubsol}
\phi_t(x,t) + G_*\!\left(x,t,\, D\phi(x,t),\, 
D^2\phi(x,t)\right) \leq 0
\end{equation}
(respectively $\geq 0$), where $G$ is given in~(\ref{eqopg}). 
A \emph{viscosity solution} of~(\ref{eqlevelset}) is a continuous 
function $v : \mathbb{R}^n \times (0,+\infty) \to \mathbb{R}$ that is simultaneously a subsolution and a supersolution.
\end{defin}
\begin{rem}
Notice that when 
$X\phi(x,t) \neq 0$, condition~(\ref{eqsubsol}) is equivalently 
written as
\begin{equation}\label{eqregsubsol}
\phi_t(x,t) + F\!\left(x,\, X\phi(x,t),\, 
X^2\phi(x,t)\right) \leq 0,
\end{equation}
as the extension $G_*$ only plays a role at the singular points where $X\phi(x,t) = 0$. At such points, the viscosity subsolution condition~(\ref{eqsubsol}) is potentially stricter 
than the condition obtained by replacing $G_*$ with $F_*$, namely
\begin{equation}\label{eqfakesubsol}
\phi_t(x,t) + F_*\!\left(x,\, X\phi(x,t),\, 
X^2\phi(x,t)\right) \leq 0.
\end{equation}
In the literature, one encounters sometimes \eqref{eqfakesubsol} to define subsolutions at characteristic points.
In the special case of the horizontal mean curvature flow operator in step two Carnot groups 
 however, when $p\,\sigma(x) = 0$, by structure condition (\ref{smatrix}) on $\sigma$ one has 
$$
G_*(x,p,A)=F_*(0,\sigma^T(x)\, A\, \sigma(x)) = \min_{|q|=1}
\bigl\{-\operatorname{Tr}\bigl((I - q \otimes q)\, \sigma^T(x)\, A\, \sigma(x)\bigr)\bigr\}.
$$
\end{rem}

We briefly recall the geometric framework of generalized flows developed by Barles and Souganidis \cite{bs} and Barles and Da Lio \cite{bdl}. In view of the structure of the equation \eqref{allen cahn eps eq}, we adapt this approach to account for singularities at characteristic points in the level-set equation for the front propagation.
This more geometric definition has been proven to be much more efficient than the standard level set approach, when dealing with singularly perturbed problems that give rise to geometric flows. 
\begin{defin}\label{generalized flow in Carnot group def} Let $F:\R^n\times\R^m\backslash\{0\}\times\mathcal S^m$ be locally bounded and satisfying (F1-2-3), and let $G$ be defined as in (\ref{eqopg}).
A family $\open$ (resp. $\close$) of open (resp. closed) subsets of $\RRn$ is called a \emph{generalized superflow} (resp. \emph{subflow}) with normal velocity $V=-F$ in the subelliptic framework if the following occurs. For any $x_0\in\RRn$, $t\in(0,T)$, $r>0$, $h>0$ so that $t+h<T$ and for any smooth function $\phi:B(x_0,r]\times[t,t+h]\rightarrow\RR$ such that (i)--(v) are satisfied
\begin{description}
  \item[(i)] $\partial_t\phi(x,s)+G^*(x,t,\nabla \phi(x,s),D^2\phi(x,s))<0$ in $B(x_0,r]\times[t,t+h]$\\ (resp. $\partial_t\phi(x,s)+G_*(x,t,\nabla \phi(x,s),D^2\phi(x,s))>0$)
  \item[(ii)] for any $s\in[t,t+h]$, $\{x\in B(x_0,r]:\phi(x,s)=0\}\neq\emptyset$ and
  $$|\nabla \phi(x,s)|\neq0\mbox{ on }\{(x,s)\in B(x_0,r]\times[t,t+h]:\phi(x,s)=0\},$$
  \item[(iii)] $\{x\in B(x_0,r]:\phi(x,t)\geq0\}\subset\Omega_t$ (resp. $\{x\in B(x_0,r]:\phi(x,t)\leq0\}\subset\mathcal{F}_t^c$),
  \item[(iv)] for all $s\in[t,t+h]$, $\{x\in \partial B(x_0,r]:\phi(x,s)\geq0\}\subset\Omega_s$ (resp. $\{x\in \partial B(x_0,r]:\phi(x,s)\leq0\}\subset\mathcal{F}_s^c$),
  \item[(v)] (only if $G$ is the mean curvature operator, see Remarks \eqref{remflow}, \eqref{remflowtwo} if there exists a pair $(x,s)\in B(x_0,r]\times[t,t+h]$ so that $\abs{X\phi(x,s)}=0$, then it holds also $\abs{X^2\phi(x,s)}=0$,
\end{description}
then $\phi$ has to satisfy
$$\{x\in B(x_0,r]:\phi(x,s)>0\}\subset\Omega_s,\quad
(\mbox{resp. }\{x\in B(x_0,r]:\phi(x,s)<0\}\subset\mathcal{F}_s^c,)$$
for every $s\in(t,t+h)$.

A family $\open$ of open subsets of $\RRn$ is called a \emph{generalized flow} with normal velocity $-F$ if $\open$ is a superflow and  $(\overline{\Omega}_t)_{t\in(0,T)}$ is a subflow.
\end{defin}
It follows immediately by Definition \ref{generalized flow in Carnot group def} that a family $\open$ of open subsets of $\RRn$ is a generalised superflow with normal velocity $-F$ if and only if $(\Omega_t^c)_{t\in(0,T)}$ is a generalised subflow with normal velocity $F$.
\begin{rem}\label{remflow}
The previous definition focuses on evolution of sets directly instead of looking at the level sets of the solutions of the partial differential equation. It does this by assuming local comparison with smooth evolutions. When we ask test functions to satisfy (i)--(iv), it is equivalent to the level set method, as a consequence of the comparison principle for the HJ equation. This is proved in \cite{bdl} if $F(x,p)=c|p|$. In the case of the mean curvature flow equation in the subelliptic setting, where we additionally require test functions to satisfy (v), we discuss it in detail in the next remark.
\end{rem}
\begin{rem}\label{remflowtwo}
Our additional condition~(v) in Definition~\ref{generalized flow in Carnot group def} is only directly needed in Section~4, specifically when handling the critical term~\eqref{estimate} via Lemma~\ref{lemtech} in the final step of that section. It is inspired by an equivalent formulation of viscosity solutions for the mean curvature flow equation, which was introduced and proved by Barles and Georgelin~\cite{bage} in the Euclidean setting, and later extended independently to step-two Carnot groups by Ferrari, Liu, and Manfredi~\cite{flm} and the authors~\cite{deso5}. The approach in~\cite{bage} has been used in the literature to prove the equivalence between generalized flows and the level-set approach in the Euclidean setting. We similarly employed it in~\cite{deso5} for the same purpose. In~\cite{deso5}, as a geometric application, we used this condition to establish in a simple way that boundaries of strictly convex sets in step-two Carnot groups become extinct in finite time under horizontal mean curvature flow, even in the presence of characteristic points. At present, it remains open to what extent this approach can be extended beyond step-two Carnot groups to more complex subelliptic geometries.
\end{rem}

\subsection{The abstract method}

We now outline the abstract method to study the asymptotic behavior of solutions to singular perturbations of semilinear reaction-diffusion equations by means of generalized sub- and superflows. We do this following the ideas explained in \cite{bs} and in \cite{bdl}.
We start recalling the concept of half relaxed limits of a locally, uniformly bounded family of functions $u^\gep:\R^n\times(0,+\infty)\to\R$. These are defined as
\begin{equation*}\begin{array}{l}
\liminf_{*\gep\to0^+}\;u^\gep(x,t):=\lim_{r\to0^+}\inf\{u^\gep(y,s):0<\gep<r,\;(y,s)\in B(x,r)\times(t-r,t+r)\}\\
\limsup_{\gep\to0^+}^*\;u^\gep(x,t):=\lim_{r\to0^+}\sup\{u^\gep(y,s):0<\gep<r,\;(y,s)\in B(x,r)\times(t-r,t+r)\}.
\end{array}\end{equation*}

Given a family $(u_\gep)_{\gep>0}$ of bounded functions in $\RRn\times[0,T]$, representing the solutions to the Cauchy problem for our reaction-diffusion equations \eqref{case1} and \eqref{case2}, we aim at proving that for any
$t\in(0,T)$, we can define two regions $\Omega_t^1$ and $\Omega_t^2$ where $u^\gep(\cdot,t)$ converges to the two stable equilibria $m_\pm$, respectively, i.e. locally uniformly
$$
u^\gep(x,t)\to
\left\{
\begin{array}{ll}
m_+&\mbox{if }(x,t)\in\Omega^1:=\bigcup_{t\in(0,T)}\Omega^1_t\times\{t\},\\
m_-&\mbox{if }(x,t)\in\Omega^2:=\bigcup_{t\in(0,T)}\Omega^2_t\times\{t\}.
\end{array}
\right.
$$
Furthermore we want to study the evolution of the interface $t\mapsto\Gamma_t$ separating $\Omega^1_t$ and $\Omega^2_t$. 
The key step is to prove that the families
$((\Omega^1_t)^c)_{t\in(0,T)}$ and $((\Omega^2_t)^c)_{t\in(0,T)}$ 
act as a superflow and a subflow, respectively, with a normal velocity determined by the problem data. 

We proceed by defining
two open sets $\Omega^1$ and $\Omega^2$, where the limit goes as expected
\begin{equation}\label{omega1-2B}
\begin{split}
\Omega^1&=\Int\Big\{(x,t)\in\RRn\times[0,T]:\liminfstard\Big[\frac{u^\gep-m_+}{\gep^\tau}\Big](x,t)\geq0\Big\}\\
\Omega^2&=\Int\Big\{(x,t)\in\RRn\times[0,T]:\limsupstard\Big[\frac{u^\gep-m_-}{\gep^\tau}\Big](x,t)\leq0\Big\},
\end{split}
\end{equation}
where $\tau\geq0$ will be chosen appropriately as $\tau=0$ in Case 1 and $\tau=1-k$ in Case 2.

For $t\in(0,T)$ we define the projections $\omegaone$ and $\omegatwo$ by
\begin{equation}\label{omega1-2t}
\begin{array}{ll}
\Omega^1_t=\{x\in\R^n:(x,t)\in\Omega^1\},\quad
\Omega^2_t=\{x\in\R^n:(x,t)\in\Omega^2\}.
\end{array}
\end{equation}
Clearly $\Omega^1=\cup_{t\in(0,T)}\Omega^1_t\times\{t\}$, $\Omega^2=\cup_{t\in(0,T)}\Omega^2_t\times\{t\}$.
Since $\Omega^1$ and $\Omega^2$ are open and disjoint subsets of $\RRn\times(0,T)$ the two step functions defined as
\begin{equation}\label{chi def}
\underline w(x,t)=\mathds{1}_{\Omega^1}-\mathds{1}_{(\Omega^1)^c},\quad
\overline w(x,t)=\mathds{1}_{(\Omega^2)^c}-\mathds{1}_{\Omega^2}
\end{equation}
are respectively lower and upper semicontinuous in $\RRn\times(0,T)$.  We observe that $\underline w$, $\overline w$ can be respectively extended by lower and upper semicontinuity to the whole of $\RRn\times[0,T)$. For simplicity of notation we still call $\underline w$ and $\overline w$ these extensions.
To analyze the asymptotics of functions $u^\epsilon$ we then follow three steps.

\noindent
1. \emph{Initialization}: we show that $\Omega^1,\Omega^2$ are not empty and we define their respective traces $\Omega^1_0$ and $\Omega^2_0$ at $t=0$ as
\begin{equation}\label{omega0}
\Omega^1_0=\{x\in\RRn:\underline w(x,0)=1\},\qquad\Omega^2_0=\{x\in\RRn:\overline w(x,0)=-1\}.
\end{equation}
\emph{2. Propagation}: we show that $\omegaone$ and $((\Omega^2_t)^c)_{t\in(0,T)}$ are respectively super and sub-flows with a normal velocity $V$ defined appropriately.

\noindent
\emph{3. Conclusion}: we conclude our asymptotics by applying Corollary \ref{cor levelset- generalized flow} to compare $(\Omega^1_t)_{t\in[0,T)}$ and $((\Omega^2_t)^c)_{t\in[0,T)}$ to the level set evolution of $(\Gamma_0,\Omega^1_0,\Omega^2_0)$, $\Gamma_0=\RRn\backslash(\Omega^1_0\cup\Omega^2_0)$.

In the next two sections we will thus consider the two cases, concentrating on the steps 1 and 2.

\section{The first asymptotic problem}

We consider the asymptotics of the solutions of the Cauchy problem for the equation (\ref{asymp eq}).
In this section, we require  $c\neq0$ if $(q(\cdot),c)$ is the solution to the traveling wave boundary value problem 
\begin{equation}\label{traveling wave equation}
\ddot q+c\dot q=f(q),\qquad \lim_{r\to\pm\infty}q(r)=m_{\pm},\qquad q(0)=m_o,
\end{equation}
which we discuss in more detail in the Appendix.
We show that suitable initial conditions generate a front moving by
 a first order normal velocity. Its level set equation is the following geometric equation ($F(p)=|p|$)
\begin{equation}\label{first ord deg eq}
v_t(x,t)+G(x,\nabla v(x,t))=v_t(x,t)+c\abs{\nabla v(x,t)\sigma(x)}=0,\qquad (x,t)\in\RRn\times(0,+\infty).
\end{equation}
Our core assumptions in this section require that $\sigma=(\sigma_{ik})_{i,k}\in C^2(\RR^n,\RR^{n\times m})\cap W^{2,\infty}(\RR^n,\RR^{n\times m})$ satisfies (\ref{A regularity}), (\ref{A non negative}), $f$ satisfies \eqref{f property} and the traveling wave $q(\cdot)$ satisfies \eqref{q property}, \eqref{q property 2}.
Notice that in our standing assumptions, there exists a positive constant $L_\sigma$ so that
$$\abs{G(x,p)-G(y,p)}\leq L_\sigma\abs{x-y}\abs{p},\qquad \abs{G(x,p)-G(x,q)}\leq L_\sigma c\abs{p-q}.$$
Then it is well known, see e.g. \cite{cil} or \cite{bcd}, that the Cauchy problem for (\ref{first ord deg eq}) is well-posed and has a unique viscosity solution for all continuous initial conditions.
Moreover the equation satisfies the comparison principle for sub and supersolutions of the equation (\ref{first ord deg eq}) (see for example \cite{cil}, \cite{bcd}).  
Therefore the standard level set approach is equivalent to the generalized front evolution defined in (\ref{generalized flow in Carnot group def}) (in this case we do not require condition (v) ).
We prove the following result for the asymptotics.

\begin{thm}\label{thm deg} Assume that the matrix map $A=\sigma(x)\sigma^{t}(x)$ and the function $f$ satisfy (\ref{A regularity}), (\ref{A non negative}) (\ref{f property}), $c\neq0$ and $q(\cdot)$ satisfies (\ref{q property}) and (\ref{q property 2}). For $\gep>0$, let $u^\gep$ be the unique solution of
\begin{equation}\label{asymp cauchy problem}
\left\{
\begin{array}{rl}
\mbox{(i)}&u^{\gep}_t(x,t)-\gep\tr(Y_\gep^2 u^\gep(x,t))+\gep^{-1}f(u^\gep)=0\quad (x,t)\in\RRn\times(0,+\infty)\\
\mbox{(ii)}&u^\gep(x,0)=g(x),\quad x\in\RRn,
\end{array}
\right.
\end{equation}
 where $g:\RRn\to[m_-,m_+]$ is a continuous function such that the sets $\Gamma_o=\{x:g(x)=m_o\}$, $\Omega_o^+=\{x:g(x)>m_0\}$, $\Omega_o^-=\{x:g(x)<m_0\}$ are nonempty subsets of $\RRn$. If the exponent $k$ in (\ref{sigmaeps def}) satisfies $0<k<1$, then
$$u^\gep(x,t)\longrightarrow
\left\{\begin{array}{rll}
m_+&\mbox{in} &\{(x,t):v(x,t)>0\},\\
m_-&\mbox{in} &\{(x,t):v(x,t)<0\},
\end{array}
\right.
$$
locally uniformly as $\gep\to0$,
where $v$ is the unique viscosity solution of
\begin{equation}\label{PCfront deg}
\left\{\begin{array}{rll}
\mbox{(i)}&v_t(x,t)+c|\nabla v(x,t)\sigma(x))|=0,\quad(x,t)\in\RRn\times(0,+\infty),\\
\mbox{(ii)}&v(x,0)=d_o(x),\quad x\in\RRn,
\end{array}\right.
\end{equation}
and $d_o$ is the signed distance to $\Gamma_o$ which is positive in $\Omega^+_o$ and negative in $\Omega_o^-$.
\end{thm}
\begin{rem}
The proofs in this and the subsequent section are necessarily technical. To maintain readability and avoid repeating well-established machinery, we outline the complete architectural structure of the arguments while relying on standard adaptations from the Euclidean literature (such as the gluing constructions of global subsolutions in Step~2 of the second lemma, and the basic initializations in the first lemma, developed in~\cite{bs, bdl, xc}). Where the proofs diverge, specifically in handling the sub-Riemannian geometry, characteristic points, and the refined Riemannian approximation bounds, i.e., the modified definitions of the local sub and super-solutions, full details are provided. Since the streamline of Sections 3 and 4 is the same, some calculations are repeated in several Lemmas. When this happens with minor modifications we will omit the full details.
\end{rem}
Before entering our proof we make some preliminary remarks and formal computations to describe some significant changes that we make into the classical argument. In \cite{bss},
for any $\gep>0$ suitable sub and supersolutions of the equation (\ref{asymp cauchy problem}-i) are constructed in a neighborhood of a level set $\{\phi(\cdot,t)=0\}$, where $\nabla\phi\neq0$, by using the formula
$\Phi^\gep(x,t)=q(\frac{z^\gep(x,t)}{\gep})$,
where $z^\gep$ is chosen appropriately.
From the equation, we obtain
\begin{equation*}
\frac{\dot q}{\gep}[z_t-\gep\tr(Y_\gep^2z)+c]+\frac{\ddot q}{\gep^2}[1-|Y_\gep z|^2]\leq0.
\end{equation*}
Since $\dot q/\gep, \;\ddot q/\gep$ are small away from the interface if $\gep$ is small, this last inequality is satisfied if in a neighborhood of the front it holds
\begin{equation}\label{heuristic1}
\begin{array}{c}
z_t(x,t)-\gep\tr(Y_\gep^2z(x,t))+c\abs{Y_\gep z(x,t)}\leq0\\
|Y_\gep z(x,t)|=1.
\end{array}\end{equation}
Let's observe the second equation in (\ref{heuristic1}) in a neighborhood of $\{\phi=0\}$. If in such neighborhood $\nabla\phi\neq0$ and by \eqref{Agep pos def},
by the method of characteristics, for $t>0$, 
$z^\epsilon(\cdot,t)$ could be the Riemannian signed distance from the front having the same sign of $\phi$. Since its horizontal gradient $X z^\gep$ might vanish even though the gradient does not, the Riemannian signed distance will still approximate the nonsmooth sub-Riemannian signed distance solving as a viscosity solution
 $$\abs{Xz}=1.$$
 Instead, we will choose to define directly
 $z^\gep=\phi\abs{Y_\gep\phi}^{-1}$. This choice will produce additional terms in the computations. Moreover 
 $$Y_\gep z^\gep=\frac{Y_\gep\phi}{\abs{Y_\gep \phi}}-\frac{\phi}{\abs{Y_\gep \phi}}\frac{Y_\gep \phi \;Y_\gep^2 \phi}{\abs{Y_\gep \phi}},
 $$
so that the second term will vanish on the front, and its denominator might also vanish in the limit as $\gep\to0$ at characteristic points, but at a rate controllable with $\gep^k$. The choice of the exponent $k$ will then become crucial. Differentiating these terms will increase the power.

We prove Theorem \ref{thm deg} by following the abstract method, as recalled in the previous section. We define two families of open sets of $\RRn$, $\omegaone$ and $\omegatwo$ by (\ref{omega1-2B}) with $\tau=0$, (\ref{omega1-2t}) and extend them at $t=0$ by ($\ref{omega0}$). 
We want to show that $\Omega_0^+=\{x:g(x)>m_o\}\subseteq\Omega_0^1$ and $\Omega_0^-=\{x:g(x)<m_o\}\subseteq\Omega_0^2$ and that they generate sub and superflows.
The proof of Theorem \ref{thm deg} is then completed by Corollary \ref{cor levelset- generalized flow} in the Appendix.
Observe that by the maximum principle $m_-\leq u^\gep\leq m_+$ for all $\gep>0$. We will thus focus on the two steps, initialization and propagation.

\emph{First step: initialization.} Here we show that $\Omega_0^+=\{x:g(x)>m_o\}\subseteq\Omega_0^1$, the other inclusion being similar.
  Let $\hat x$ so that $g(\hat x)>m_o$ and find $r,\sigma>0$ so that $g(x)\geq m_o+\sigma$  for any $x\in B(\hat x, r)$.
 This means that
\begin{equation}\label{deg lower estimate uepsilona}
u^\gep(x,0)=g(x)\geq (m_o+\sigma)\mathds1_{B(\hat x,r)}(x)+m_-\mathds1_{B(\hat x,r)^c}(x).
\end{equation}
Now we introduce the function $\Phi:\RRn\times[0,T]\to\RR$ defined as
\begin{equation}\label{deg Phi def}
\Phi(x,t)=r^2-|x-\hat x|^2-Ct,
\end{equation}
with $C>0$ a suitable constant that will be chosen later. 
The main part of the proof that $\hat x\in\Omega^1_0$ is contained in the following two lemmas. The first one follows \cite{bdl}, see also the original argument of Chen \cite{xce,xc} and we sketch it below for the reader's convenience. 

\begin{lem}\label{deg initialization first lemma} Under the assumption of Theorem \ref{thm deg} we have that for any $\beta>0$ there exist $\tau=\tau(\beta)>0$ and $\bar\gep=\bar\gep(\beta)$ such that, for all $0<\gep\leq\bar\gep$, we have
$$u^\gep(x,t_\gep)\geq(m_+-\beta)\mathds1_{\{\Phi(\cdot,0)\geq\beta\}}(x)+m_-\mathds1_{\{\Phi(\cdot,0)<\beta\}}(x),\quad x\in\RRn,$$
where $t_\gep=\tau\gep$.
\end{lem}
\begin{proof}
The idea is that for $\gep$ small enough and for short time the diffusion term in (\ref{asymp cauchy problem}-i) is negligible.
Let $\chi=\chi(\tau,\xi)\in C^2([0,+\infty)\times\RR)$ be the solution of
\begin{equation*}
\left\{
\begin{array}{l}
\dot\chi(\tau,\xi)+f(\chi(\tau,\xi))=0,\quad \tau>0,\\
\chi(0,\xi)=\xi.
\end{array}
\right.
\end{equation*}
By (\ref{f property}) and comparison of ODEs, $\chi$ satisfies 
\begin{equation}\tag{$\chi1$}\label{property one chi}
\chi_\xi(\tau,\xi)>0,\quad\mbox{ in }[0,+\infty)\times\RR
\end{equation}
and moreover there exists $\tau_o=\tau_o(\beta)>0$ such that, for all $\tau\geq\tau_o$
\begin{equation}\tag{$\chi2$}\label{property two chi}
\begin{array}{l}
\chi(\tau,\xi)\geq m_+-\beta,\quad\forall\,\xi\geq m_o+\frac\sigma2.
\end{array}
\end{equation}
Finally, for any $C>\max\{|m_--m_0|,|m_+-m_0|\}$, we have that $\chi(\tau,\xi)\in[m_0-C,m_0+C]$ for all $\xi\in[-C,C]$, $\tau\geq0$. Therefore there exists a constant $M_{C,\tau}>0$ such that, for any $\xi\in[-C,C]$ and $\gep$ small enough,
\begin{equation}\tag{$\chi3$}\label{property three chi}
|\chi_{\xi\xi}(\tau,\xi)|\leq M_{C,\tau}\chi_\xi(\tau,\xi).
\end{equation}
Now define the function $\underline u^\gep$ in $\RRn\times[0,T]$ as
$$\underline u^\gep(x,t)=\chi\Big(\frac{t}{\gep},\psi(\Phi(x,0))-Kt\Big),$$
where $\Phi$ is defined in (\ref{deg Phi def}), $K$ is a constant to be chosen and $\psi$ is a nondecreasing smooth function in $\RR$ such that
$$\psi(z)=\left\{\begin{array}{ll}
m_-&\mbox{if }z\leq 0,\\
m_o+\sigma&\mbox{if }z\geq\beta\wedge\frac{\sigma}{2}.
\end{array}\right.
$$
For a suitable choice of the constant $K$, and $\gep$ small enough $\underline u^\gep$ is a subsolution of (\ref{asymp cauchy problem}-i) in $\RRn\times(0,\tau_o\gep)$. Indeed, by (\ref{property one chi}) and (\ref{property three chi}) we obtain
$$\begin{array}{ll}
\underline u_t^\gep-\gep\tr(Y_\gep^2\underline u^\gep)+\frac{f(\underline u^\gep)}{\gep}&=-K\chi_\xi-\gep\big[(\chi_{\xi\xi}(\psi')^2+\chi_\xi\psi'')\abs{Y_\gep\Phi}^2+\chi_\xi\psi' \tr(A^\gep(x)D^2\Phi)\big]\\
&\quad-\gep\chi_\xi\psi'\tr((\nabla \Phi D\sigma(x))\sigma(x))\\
&\leq\chi_\xi(-K+C_0\gep)\leq0
\end{array}$$
for $\gep$ small enough.
Moreover, since $\underline u^\gep(\cdot,0)=\psi(\Phi(\cdot,0))\leq u^\gep(\cdot,0)$ we can apply the maximum principle in order to obtain $\underline u^\gep(x,t)\leq u^\gep(x,t)$ for any $(x,t)\in\RRn\times[0,\tau_o\gep]$. In particular
$$\begin{array}{ll}
u^\gep(x,\tau_o\gep)&\geq\underline u^\gep(x,\tau_o\gep)=\chi(\tau_o,\psi(\Phi(x,0))-K\tau_o\gep)\\
                    &=\chi(\tau_o,m_o+\sigma-K\tau_o\gep)\qquad\mbox{ if }\Phi(x,0)\geq\beta.
\end{array}$$
From this last inequality and property (\ref{property two chi}) of $\chi$ we get $u^\gep(x,\tau_o\gep)\geq m_+-\beta$, for $\gep$ small enough.
\end{proof}
Our first main result of the section is the following Lemma. The most delicate new parts are in steps 1 and 3 of the proof.
\begin{lem}\label{deg initialization second lemma}There exist $\bar h=\bar h(r,\hat x),\,\bar\beta=\bar\beta(r,\hat x)>0$ independent of $\gep$ such that if $\beta\leq\bar\beta$ and $\gep\leq\bar\gep(\beta)\wedge1$, then there exists a subsolution $\omega^{\gep,\beta}$ of (\ref{asymp cauchy problem}-i) in $\RRn\times(0,\bar h)$ that satisfies
$$\omega^{\gep,\beta}(x,0)\leq(m_+-\beta)\mathds1_{\{\Phi(\cdot,0)\geq\beta\}}(x)+m_-\mathds1_{\{\Phi(\cdot,0)<\beta\}}(x),\quad x\in\RRn.$$
If  moreover $(x,t)\in B(\hat x,r)\times(0,\bar h)$ and $\Phi(x,t)>3\beta$, then $$\liminfstard\omega^{\gep,\beta}(x,t)\geq m_+-2\beta.$$
\end{lem}
\begin{proof}
Let $(q,c)$ be the pair solution of the traveling wave equation (\ref{traveling wave equation}). 
Let $\Phi$ be the smooth function defined in (\ref{deg Phi def}) with $r>0$ assigned at the beginning of the initialization step. We pick $\gamma,C,\bar h$ so that $\gamma\leq r^2/3$, $C\geq1$
\begin{equation}\label{C assumption deg}
C\geq(c+1)\sup_{Q_{\gamma,r^2/3}} \big(|\nabla \Phi\sigma|^2+|\nabla\Phi|^2\big)^{1/2}>0.
\end{equation}
and $ \bar h\leq r^2/(3C)$, where $Q_{\gamma,s}=\{(x,t):|\Phi(x,t)|\leq\gamma,\,0\leq t\leq s\}$.
By definition, we have that
$$|Y_\gep\Phi(x,t)|\geq \epsilon^k|\nabla\Phi(x,t)|\neq0\qquad\mbox{for any }(x,t)\in Q_{\gamma,\bar h},\,\gep>0.$$
Adapting the method in \cite{bs} and \cite{bdl} we now construct a subsolution of (\ref{asymp cauchy problem}-i) by steps.

1. Construction of a strict subsolution of (\ref{asymp cauchy problem}-i) in $Q_{\gamma,\bar h}$.\\ 
For a small $\beta, \gep>0$, we define a smooth function $v^\gep$ in $Q_{\gamma,\bar h}$, our candidate subsolution, as
\begin{equation}\label{eqtest3a}
v^{\gep}(x,t)=q\Big(\frac{\Phi(x,t)-2\beta}{\gep|Y_\epsilon \Phi|}\Big)-2\beta.\end{equation}
Putting $v^\gep$ inside (\ref{asymp cauchy problem}-i) we obtain, here  we used the notation $s=\frac{\Phi-2\beta}{\gep}$ for convenience,
$$\begin{array}{ll}
v_t^{\gep}-\gep\tr(Y^2_\gep v^\gep)+\frac{f(v^\gep)}{\gep}=\frac{\dot q\Phi_t}{\gep|Y_\epsilon \Phi|}-\frac{\ddot q}{\gep}-\frac{\dot q}{|Y_\gep \Phi|}\tr(Y^2_\gep \Phi)+\frac{f(q-2\beta)}{\gep}\\
\qquad-\gep\dot qs\tr\big(Y_\gep^2(\frac{1}{\abs{Y_\gep\Phi}})\big)-\gep\ddot qs^2\big|Y_\gep(\abs{Y_\gep\Phi}^{-1})\big|^2
-2\Big(\frac{\ddot qs}{\abs{Y_\gep\Phi}}+\dot q\Big)\langle Y_\gep\left(\frac{1}{\abs{Y_\gep\Phi}}\right),Y_\gep\Phi\rangle\\
\qquad\leq\frac{I_\gep}{\gep}+\frac{I\!I_\gep}{\gep}+I\!I\!I_\gep
\end{array}$$
where we use a Taylor expansion of $f$ around $q(r)$ and the traveling wave equation so that
\begin{equation}\label{I-II-III}\begin{array}{c}
I_\gep=-\ddot q+f(q)=c\dot q(\frac{s}{\abs{Y_\gep\Phi}}) \\
I\!I_\gep=\frac{\dot q\Phi_t}{|Y_\epsilon \Phi|}-\gep\frac{\dot q}{|Y_\gep \Phi|}\tr(Y^2_\gep \Phi)-2\beta f'(q)+2\norm{f''_{|[m_-,m_+]}}_\infty\beta^2\\
{I\!I\!I}_\gep=-\gep\dot qs\tr\big(Y_\gep^2(\frac{1}{\abs{Y_\gep\Phi}})\big)-\gep\ddot qs^2\big|Y_\gep(\abs{Y_\gep\Phi}^{-1})\big|^2
-2\Big(\frac{\ddot qs}{\abs{Y_\gep\Phi}}+\dot q\Big)
\frac{\langle Y_\gep\Phi Y^2_\gep\Phi,Y_\gep\Phi\rangle}{\abs{Y_\gep\Phi}^3}.
\end{array}
\end{equation}
Observe that
$$\begin{array}{ll}
I\!I_\gep&=\frac{\dot q}{|Y_\gep\Phi|}[-C-2\gep\tr(A^\gep(x)+\left((x-\hat x)^TD\sigma(x)\right)\sigma(x))]-2\beta f'(q)+2\norm{f''_{|[m_-,m_+]}}_\infty\beta^2 \\
&\leq\frac{\dot q}{|Y_\gep\Phi|}(-C+C_1\gep)-2\beta f'(q)+2\norm{f''_{|[m_-,m_+]}}_\infty\beta^2
\end{array}$$
and
$$
I\!I\!I_\gep\leq-\gep\dot qs\tr\big(Y_\gep^2(\frac{1}{\abs{Y_\gep\Phi}})\big)+\gep\frac{\abs{\ddot q}s^2}{\abs{Y_\gep\Phi}^4}\norm{Y^2_\gep\Phi}^2+2\Big(\frac{\abs{\ddot q}\abs{s}}{\abs{Y_\gep\Phi}}+\dot q\Big)
\frac{\norm{ Y_\gep^2\Phi}}{\abs{Y_\gep\Phi}}.
$$
By assumption (\ref{q property}) on $q(\cdot)$, and computing third derivatives
$$\begin{array}{ll}
\gep I\!I\!I_\gep&\leq-\gep^2\dot qs\tr\big(Y_\gep(Y_\gep\Big(\frac{1}{\abs{Y_\gep\Phi}}\Big)\big)+C_1(\frac{\gep^2}{\abs{Y_\gep\Phi}^2}+\frac{\gep}{\abs{Y_\gep\Phi}})\\
&\leq-\gep^2\dot qs\tr\big(Y_\gep\Big(\frac{-Y_\gep\Phi Y_\gep^2\Phi}{\abs{Y_\gep\Phi}^3}\Big)\big)+C_2(\gep^{2-2k}+\gep^{1-k})\\
&\leq C_3\left(\gep^2\dot q\frac{|s|}{|Y_\gep\Phi|^3}+\gep^{1-k}\right)\leq C_4\left(\gep^{2-2k}+\gep^{1-k}\right)=o_\gep(1).
\end{array}$$
Combining all the estimates above and using the assumption on $C$ in (\ref{C assumption deg}) we get
\begin{equation*}
\begin{array}{ll}
v_t^{\gep}-\gep\tr(Y_\gep^2 v^\gep)+\frac{f(v^\gep)}{\gep}\!\!&\leq\frac{1}{\gep}\Big[\dot q\Big(c\!-\frac{C}{|Y_\gep\Phi|}+o_\gep(1)\Big)\!-2\beta f'(q)\!+\!2\beta^2\norm{f''_{|[m_-,m_+]}}_\infty\!\!+o_\gep(1)\Big]\\
&\leq\frac{1}{\gep}\Big[\dot q(-1+o_\gep(1))-2\beta f'(q)+2\beta^2\norm{f''_{|[m_-,m_+]}}_\infty+o_\gep(1)\Big],\end{array}\end{equation*}
and thus, for $\gep$ small enough,
\begin{equation}\label{computation1}
v_t^{\gep}-\gep\tr(Y_\gep^2 v^\gep)+\frac{f(v^\gep)}{\gep}\leq\frac{1}{\gep}\Big[-\frac{\dot q}{2}-2\beta f'(q)+2\beta^2\norm{f''_{|[m_-,m_+]}}_\infty+o_\gep(1)\Big].
\end{equation}
To prove that $v^\gep$ is a subsolution of (\ref{asymp cauchy problem}-i) it remains to see that the right hand side of (\ref{computation1}) is non positive. To do this we recall that $f'(m_\pm)\geq d>0$ and $q(r)\to m_\pm$ if $r\to\pm\infty$. This means that there exists an $M>0$ (independent of $\beta$) such that
$$f'(q(r))\geq\frac{d}{2},\quad\mbox{for any }|r|\geq M.$$
Moreover we can choose $\beta$ small enough in order to get
$$4\beta\norm{f''_{|[m_-,m_+]}}_\infty\leq d.$$
Therefore from \eqref{computation1}, if $\frac{|\Phi-2\beta|}{\gep|Y_\gep\Phi|}\geq M$, we have that ($\dot q\geq0$)
$$v_t^{\gep}-\gep\tr(Y_\gep^2 v^\gep)+\frac{f(v^\gep)}{\gep}\leq
\frac{1}{\gep}\Big[-\frac{\beta d}{2}+o_\gep(1)\Big]\leq-\frac{\beta d}{4\gep}<0$$
for $\gep$ small enough. If instead $\frac{|\Phi-2\beta|}{\gep|Y_\gep\Phi|}< M$ and we choose $K>0$ (which depends on $M$) so that  $\dot q(r)\geq K$ for any $|r|\leq M$, we can find $\mu>0$ so that
 \begin{equation*}\begin{array}{ll}
 v_t^{\gep}-\gep\tr(Y_\gep^2 v^\gep)+\frac{f(v^\gep)}{\gep}&\leq\frac{1}{\gep}[-\frac{K}{2}+2\beta\norm{f'}_\infty+2\beta^2\norm{f''}_\infty+o_\gep(1)]\leq-\frac{\mu}{\gep}<0.\end{array}\end{equation*}
 for $\beta$ small compared with $K$, and $\gep$ small enough.

 2. Construction of a global subsolution of (\ref{asymp cauchy problem}-i).
 
 \noindent
 This part follows closely Lemma 4.4 in \cite{bs} and we only insist on the main adaptations for the reader's convenience.
 Once we have proved that $v^\gep$ is a strict subsolution of (\ref{asymp cauchy problem}-i) in $Q_{\gamma,\bar h}$ we define, for each $(x,t)\in\{(x,t)\in\RRn\times[0,\bar h]:\Phi(x,t)\leq\gamma\}$,
$$\bar v^\gep(x,t)=
\left\{\begin{array}{lll}
\max(v^\gep(x,t),m_-)&\mbox{if}&-\gamma\leq\Phi(x,t)\leq\gamma,\\
m_-&\mbox{if}&\Phi(x,t)<-\gamma.
\end{array}
\right.$$
We obtain that $\bar v^\gep$ is a continuous viscosity subsolution of (\ref{asymp cauchy problem}-i) in $\{(x,t)\in\RRn\times[0,\bar h]:\Phi(x,t)\leq\gamma\}$, for $\gep$ sufficiently small.  On $\{(x,t):\abs{\Phi(x,t)}<\gamma\}$, $\bar v^\gep$ is the supremum of two subsolutions so this is well known. 
We can moreover show that $\bar v^\gep(x,t)=m_-$ 
if $(x,t)$ is such that $-\gamma\leq\Phi(x,t)\leq-\gamma/2$, 
so the statement becomes obvious. Indeed by properties (\ref{q property}) of $q(\cdot)$, we have that, for some $\rho>0$ and $\gep$ sufficiently small,
$$v^\gep(x,t)\leq q\Big(-\frac{\gamma+4\beta}{2\gep\abs{Y_\gep\Phi}}\Big)-2\beta\leq m_-+a e^{-\frac{b(\gamma+4\beta)}{2\gep\rho}}-2\beta\leq m_-.$$
Therefore $\bar v^\gep$ is a subsolution of (\ref{asymp cauchy problem}-i) in $\{(x,t)\in\RRn\times[0,\bar h]:\Phi(x,t)\leq\gamma\}$.

Finally define the function $\omega^{\gep,\beta}:\RRn\times[0,\bar h]\to\RR$ as
$$
\omega^{\gep,\beta}(x,t)=
\left\{\begin{array}{ll}
\psi(\Phi(x,t))\bar v^{\gep}(x,t)+(1-\psi(\Phi(x,t)))(m_+-\beta)&\mbox{if }\Phi(x,t)<\gamma,\\
m_+-\beta&\mbox{if }\Phi(x,t)\geq\gamma,
\end{array}
\right.
$$
where $\psi:\RR\to\RR$ is a smooth function such that $\psi'\leq0$ in $\RR$, $\psi=1$ in $(-\infty,\gamma/2]$, $0<\psi<1$ in $(\gamma/2,3\gamma/4)$ and $\psi=0$ in $[3\gamma/4,+\infty)$. We show that $\omega^{\gep,\beta}$ is a subsolution of (\ref{asymp cauchy problem}-i) in $\RRn\times[0,\bar h]$ for $\gep$ small. The only 
check needs to be made in
$\{(x,t)\in\RRn\times(0,\bar h):\gamma/2\leq \Phi(x,t)\leq3\gamma/4\}$. 
If we take $\beta$ so that  $2\beta<\gamma/4$, for some $\rho>0$,
 we have that 
 $$\Phi(x,t)-2\beta>\gamma/4,\qquad v^\gep(x,t)\geq  m_+-a e^{-\frac{b\gamma}{4\gep\rho}}-2\beta$$
  therefore $\bar v^\gep(x,t)=v^\gep(x,t)$, for $\beta$ and then $\gep$ sufficiently small.
  Checking the equation we obtain
\begin{multline}\label{computation2}
\omegaeb_t-\gep\tr(Y_\gep^2\omegaeb)+\frac{f(\omegaeb)}{\gep}\leq [\psi'\;(\Phi_t-\gep\tr(Y^2_\gep\Phi))-\gep\psi''|Y_\gep\Phi|^2](v^\gep-(m_+-\beta))\\ +\psi\;(v^\gep_t-\gep\tr(Y^2_\gep v^\gep))-2\gep\psi'\;\langle Y_\gep \Phi,Y_\gep v^\gep\rangle+\frac{f(\omegaeb)}{\gep}.
\end{multline}
Since by definition, $v^\gep-(m_+-\beta)\leq-\beta$, $f$ is convex in a neighborhood of $m_+$ and $v^\gep$ is a strict subsolution of (\ref{asymp cauchy problem}-i) in $Q_{\gamma,\bar h}$, the inequality in (\ref{computation2}) becomes
\begin{equation*}
\begin{split}
\omegaeb_t-\gep\tr(Y_\gep^2\omegaeb)+\frac{f^\gep(\omegaeb)}{\gep}&\leq-C\psi'\;(v^\gep-(m_+-\beta))+\psi\;\Big(v^\gep_t-\gep\tr(Y^2_\gep v^\gep)+\frac{f(v^\gep)}{\gep}\Big)\\
&\qquad\qquad\qquad\qquad\qquad\qquad+C_3\gep+(1-\psi)\frac{f(m_+-\beta)}{\gep}\\
&\leq-\frac{\psi\bar\mu(\beta)}{\gep}+C_3\gep+(1-\psi)\frac{f(m_+-\beta)}{\gep}\\
&=-\frac{1}{\gep}\Big(\psi\bar\mu(\beta)+(1-\psi)(-f(m_+-\beta))\Big)+C_3\gep\\
&=-\frac{\hat\alpha_\beta}{\gep}+C_3\gep\leq0,\quad\mbox{for $\beta$ and then $\gep$ small enough}.
\end{split}\end{equation*}
In the last chain of inequalities $\bar\mu(\beta)=\mu\wedge\beta d/4$ and in the last line we have used the fact that $\hat \alpha_\beta$ is strictly positive as it is a convex combination of strictly positive terms.

3. We end proving the estimates for $\omega^{\gep,\beta}(x,0)$, $x\in\RRn$. Observe that for a given $\beta$, if $(x,t)\in Q_{\gamma,\bar h}$ satisfies $\Phi(x,t)<\beta$, then for some $\rho>0$ and the properties of the traveling wave
$$v^\gep(x,t)\leq q\Big(\frac{-\beta}{\gep\abs{Y_\gep\Phi(x,t)}}\Big)-2\beta\leq m_-+ae^{-\frac{b\beta}{\gep\rho\Phi}}-2\beta\leq m_-
$$
 for $\gep$ small enough. This means that for any $(x,t)\in Q_{\gamma,\bar h}$ it holds
 $$v^\gep(x,t)\leq (m_+-\beta)\mathds1_{\{\Phi\geq\beta\}}(x,t)+m_-\mathds1_{\{\Phi<\beta\}}(x,t),$$
 for $\gep$ small enough. By definition  
 this inequality still holds for $\bar v^\gep$ and $\omegaeb$ in their domain of definition. 
Finally we just remark that, with a similar reasoning we can prove that if $(x,t)\in B(\hat x,r)\times(0,\bar h)$ and $\gamma\geq\Phi(x,t)>3\beta$, then
$$v^\gep(x,t)\geq q\left(\frac{\beta}{\gep\abs{Y_\gep\Phi}}\right)-2\beta\geq m_+-ae^{-\frac{b\beta}{\gep\rho\Phi}}-2\beta.$$
Hence $\liminfstard\omegaeb(x,t)\geq m_+-2\beta$ for any $(x,t)\in Q_{\gamma,\bar h}$.
\end{proof}
We now conclude the initialization step with the classical argument combining the two Lemmas. For a small $\beta\leq\bar\beta$ and $\gep\leq\bar\gep(\beta)\wedge1$, we proved the existence of a subsolution $\omega^{\gep,\beta}$ of (\ref{asymp cauchy problem}-i) in $\RRn\times(0,\bar h)$ such that
$$\omega^{\gep,\beta}(x,0)\leq u^\gep(x,t_\gep),\qquad x\in\RRn.$$
Thus by the maximum principle
$$\omega^{\gep,\beta}(x,s)\leq u^\gep(x,\tau\gep+s),\qquad (x,s)\in\RRn\times [0,\bar h].$$
By taking the limit and applying Lemma \ref{deg initialization second lemma} we obtain
$$\liminfstard u^\gep(x,s)\geq m_+-2\beta,\quad\hbox{if }\Phi(x,s)>3\beta.$$
As $\beta$ is arbitrary and does not depend on $\bar h$ we end up with $x\in\Omega^1_s$, for $s\in[0,\bar h]$, if $\Phi(x,s)>0$. 
Therefore, for $\hat x\in\Omega^+_0$, we can find $\hat h,\eta>0$ so that, by definition of $\Phi$, $B(\hat x,\eta)\subset \Omega^1_s$ for all $s\in [0,\hat h]$. This also implies $\hat x\in \Omega^1_0$ which is what we want. 

\emph{Second step: propagation.} In this step we want to show that $\omegaone$ and $((\Omega^2_t)^c)_{t\in(0,T)}$ are respectively super and subflow with normal velocity $V$ determined by
$$G(x,p)=c|p\sigma(x)|.$$
 Since the two proofs are similar we'll do only the one for $\omegaone$.
Let $x_0\in\RRn$, $t\in(0,T)$, $r>0$, $h>0$ so that $t+h<T$. Suppose that $\phi:\RRn\times[0,T]\to\RR$ is a smooth function 
satisfying (i)-(iv) of Definition \ref{flow-viscositysol}, and let $\tilde C>0$ be such that
$$\phi_t(x,s)+c|\nabla \phi(x,s)\sigma(x)|\leq-\tilde C, \hbox{for all }(x,s)\in B(x_0,r]\times[t,t+h].$$
We have to show that for every $s\in(t,t+h)$,
$$\{x\in B(x_0,r]:\phi(x,s)>0\}\subset\Omega_s^1.$$
Using the assumptions and the definition of $(\Omega^1_t)_{t\in(0,T)}$ this is equivalent to prove that for any $x\in B(x_0,r)$, $s\in(t,t+h)$ such that $\phi(x,s)>0$, we have
\begin{equation}\label{propagation thesis}\liminfstard u^\gep(y,\tau)\geq m_+\end{equation}
for $(y,\tau)$ in a neighborhood of $(x,s)$ where, for any $\gep>0$, $u^\gep$ is the solution of the Cauchy problem (\ref{asymp cauchy problem}).  The proof proceeds like the second Lemma of the first step with the difference that here we have to construct a subsolution of (\ref{asymp cauchy problem}-i) only in the cylinder $B(x_0,r)\times (t,t+h)$.
 It is enough to prove the following lemma which plays the role of Lemma \ref{deg initialization second lemma}.
 
\begin{lem}\label{deg lemma propagation} Let $\phi$ be a smooth function, satisfying (i)-(iv) of Definition \ref{flow-viscositysol} as above.
There exists $\bar\beta$ small enough such that, for all $\beta\leq\bar\beta$ and $\gep\leq\bar\gep(\beta)$ then there is a viscosity subsolution $\omega^{\gep,\beta}$ of (\ref{asymp cauchy problem}-i) in $B(x_0,r)\times(t,t+h)$ that satisfies,

 1. $\omega^{\gep,\beta}(x,t)\leq(m_+-\beta)\mathds1_{\{\phi(\cdot,t)\geq\beta\}}(x)+m_-\mathds1_{\{\phi(\cdot,t)<\beta\}}(x),\quad\mbox{for all } x\in B(x_0,r],$

2. $\omega^{\gep,\beta}(x,s)\leq(m_+-\beta)\mathds1_{\{\phi(\cdot,s)\geq\beta\}}(x)+m_-\mathds1_{\{\phi(\cdot,s)<\beta\}}(x),\quad\mbox{for all } x\in \partial B(x_0,r]$, $s\in[t,t+h]$,

3. if $(x,s)\in B(x_0,r]\times[t,t+h]$ satisfies $\phi(x,s)>3\beta$, then
$$\liminfstard\omega^{\gep,\beta}(x,s)\geq m_+-2\beta.$$
\end{lem}
 \begin{proof} 
 The function $\phi$ in the statement is of course not explicit but we only need it locally, and it satisfies conditions (i)-(iv) of Definition \ref{generalized flow in Carnot group def}.
 Observe that by property (ii), there exists $\gamma>0$ such that $|\nabla\phi(x,s)|\neq0$ in the set $Q_\gamma=\{(x,s)\in B(x_0,r]\times[t,t+h]:|\phi(x,s)|\leq\gamma\}$. As in Remark \ref{remgradest} $|Y_\gep\phi(x,s)|\geq\gep^k|\nabla\phi(x,s)|\neq0$ for any $(x,s)\in Q_\gamma$, $\gep>0$. 
 We define first a function $v^{\gep}$ in $Q_\gamma$ as
\begin{equation}\label{eqtest3b}
v^{\gep}(x,s)=q\Big(\frac{\phi(x,s)-2\beta}{\gep\abs{Y_\gep\phi(x)}}\Big)-2\beta.\end{equation}
By the same computations as in step 1 of Lemma \ref{deg initialization second lemma},
$$v^\gep_t(x,s)-\gep\tr(Y_\gep^2v^\gep(x,s))+\frac{f(v^\gep(x,s))}{\gep}\leq\frac{I_\gep}{\gep}+\frac{I\!I_\gep}{\gep}+{I\!I\!I}_\gep,\quad (x,s)\in Q_\gamma,$$
where $I_\gep$, $I\!I_\gep$, ${I\!I\!I}_\gep$ are the same terms as in (\ref{I-II-III}) with $\Phi$ replaced by $\phi$. 
Only ${I\!I\!I}_\gep$ has the extra term (it was zero before)
$$\dot q s\partial_t(|Y_\gep\phi|^{-1}|)=\dot q s\Big(-\frac{\langle Y_\gep\phi, Y_\gep\partial_t\phi\rangle}{\abs{Y_\gep\phi}^3}\Big)\leq \frac{C_1}{\abs{Y_\gep\phi}},$$
where we put again $s=(\phi-2\beta)/\gep$ and use \eqref{q property}.
Again as in Lemma \ref{deg initialization second lemma} 
$$I_\gep= c\dot q$$
and by property (i) of $\phi$
 $$I\!I_\gep\leq\frac{\dot q}{|Y_\gep\phi|}[-\tilde C-c|X\phi|-\gep\tr(Y_\gep^2\phi)]-2\beta f'(q)+2\norm{f''}\beta^2.$$
As for the terms in $I\!I\!I_\gep$, also in this case
$$\gep I\!I\!I_\gep\leq 2C_4\gep^{1-k}+\gep \frac{C_1}{\abs{Y_\gep\phi}}
\leq C_5\gep^{1-k}=o_\gep(1).$$
Since 
$$\abs{Y_\gep\phi}-\abs{X\phi}=\frac{\gep^{2k}|\nabla\phi|}{\abs{Y_\gep \phi}+\abs{X\phi}}\leq C_6\gep^{k},$$
computing the equation we get
\begin{equation}\label{computation6}\begin{array}{ll}
v_t^{\gep}-\gep\tr(Y_\gep^2 v^\gep)+\frac{f(v^\gep)}{\gep}&\leq\frac{1}{\gep}\Big[\frac{\dot q}{|Y_\gep\phi|}\Big(-\tilde C+c(|Y_\gep\phi|-|X\phi|)-\gep\tr(Y_\gep^2\phi)\Big)\\
&\quad-2\beta f'(q)+2\norm{f''}_\infty\beta^2+o_\gep(1)\Big]\\
&\leq\frac{1}{\gep}\Big[-\frac{\dot q\tilde C}{2\abs{Y_\gep\phi}}-2\beta f'(q)+2\beta^2\norm{f''}_\infty+o_\gep(1)\Big]
\end{array}\end{equation}
for $\gep$ small enough. To prove that the right hand side of (\ref{computation6}) is strictly negative we proceed similarly as from \eqref{computation1} in the proof of Lemma \ref{deg initialization second lemma}. 
Again we may suppose that
$$\begin{array}{ll}
f'(q(r))\geq\frac{d}{2},\;\quad\mbox{for }\abs{r}\geq M,\quad
\dot q(r)\geq K>0,\mbox{for }\abs{r}\leq M,\\
4\beta\norm{f''}\leq d,\;\hbox{for }\beta \hbox{ small enough}.
\end{array}$$
 Thus, if $\abs{\phi(x,s)-2\beta}\geq M\gep\abs{Y_\gep\phi(x,s)}$ 
 we have ($\dot q\geq0$)
$$v_t^{\gep}-\gep\tr(Y_\gep^2 v^\gep)+\frac{f(v^\gep)}{\gep}
\leq\frac{1}{\gep}\big(-\frac{\beta d}{2}+o_\gep(1)\big)\leq -\frac{\beta d}{4\gep}$$
 for $\gep$ small enough. If instead  $\abs{\phi(x,s)-2\beta}<M\gep\abs{Y_\gep\phi}$, 
 we assume that $\beta$ satisfies
$$2\beta\big(\norm{f'}_\infty+\beta\norm{f''}_\infty+1\big)\sup_{(x,s)\in Q_\gamma}[\abs{X\phi(x,s)}+\abs{\nabla\phi(x,s)}]\leq\frac{K\tilde C}{2}.$$
Then the inequality in (\ref{computation6}) becomes
 $$\begin{array}{ll}
 v_t^{\gep}-\gep\tr(Y_\gep^2 v^\gep)+\frac{f(v^\gep)}{\gep}&\leq\frac{1}{\gep}\Big[-\frac{K\tilde C}{2\abs{Y_\gep\phi}}+2\beta\norm{f'}_\infty+2\beta^2\norm{f''}_\infty+o_\gep(1)\Big]\\
 &\leq\frac{1}{\gep}(-2\beta+o_\gep(1))\leq-\frac{\beta}{\gep},
 \end{array}$$
 for $\gep$ small enough compared with $\beta$ fixed.
 Now that we have proved $v^\gep$ is a strict subsolution of (\ref{asymp cauchy problem}-i) in $Q_\gamma$, the extension of $v^\gep$ to a global subsolution $\omega^{\gep,\beta}$ in $B(x_o,r]\times[t,t+h]$ and the proof that such a function satisfies the conclusions 1, 2, 3, follows more closely \cite{bs} and with basically no changes steps 2 and 3 of Lemma \ref{deg initialization second lemma} by exchanging $\Phi$ with $\phi$. We therefore omit the details.
\end{proof}
\begin{proof}[Proof of Theorem \ref{thm deg}]
We need to prove (\ref{propagation thesis}). Given $\beta<\hat\beta$, if $x\in B(x_o,r]$ and $\phi(x,t)\geq\beta>0$ then, by property (iii) of $\phi$, $x\in\Omega^1_t$, i.e.
$$\liminfstard u^\gep(x,t)> m_+-\beta.$$
 This means that there exists an $\gep_{x,t}=\gep_{x,t}(\beta)>0$ such that, for all $\gep\leq\gep_{x,t}$ and $(y,\tau)\in B(x,\gep_{x,t})\times(t-\gep_{x,t},t+\gep_{x,t})$, we have $u^\gep(y,\tau)\geq m_+-\beta$. By compactness of $\{x\in B(x_o,r]:\phi\geq\beta\}$ we can select an $\bar\gep>0$, depending only on $\beta$, so that, for all $\gep\leq\bar\gep$, and $x\in\{z\in B(x_0,r]:\phi(z,t)\geq\beta\}$ we have $u^\gep(x,t)\geq m_+-\beta$. Therefore
$$u^\gep(x,t)\geq(m_+-\beta)\mathds1_{\{\phi(\cdot,t)\geq\beta\}}(x)+m_-\mathds1_{\{\phi(\cdot,t)<\beta\}}(x).$$
for all $\gep\leq\bar\gep(\beta)$, $x\in B(x_0,r]$. By using (iv) of $\phi$ in a similar way, we can also obtain that, for $\gep$ small enough,
$$u^\gep(x,s)\geq(m_+-\beta)\mathds1_{\{\phi(\cdot,s)\geq\beta\}}(x)+m_-\mathds1_{\{\phi(\cdot,s)<\beta\}}(x),$$
for any $(x,s)\in \partial B(x_0,r]\times[t,t+h]$.
Combining these inequalities with those in 1. and 2. in the statement of Lemma \ref{deg lemma propagation} we can conclude, by the parabolic maximum principle, that
$$\omega^{\gep,\beta}(x,s)\leq u^\gep(x,s),\quad\mbox{for all }(x,s)\in B(x_0,r]\times[t,t+h].$$
Moreover by property 3. of $\omegaeb$ in Lemma \ref{deg lemma propagation},
$$\liminfstard u^\gep(x,s)\geq m_+-2\beta$$
for every $(x,s)\in B(x_0,r]\times[t,t+h]$ such that $\phi(x,s)>3\beta$. Since $\beta$ is arbitrary we can now send $\beta$ to zero in order to obtain that
$$\liminfstard u^\gep(x,s)\geq m_+,\qquad \hbox{for all }(x,s)\in B(x_0,r)\times(t,t+h) \hbox{ and }\phi(x,s)>0.$$
  Thus (\ref{propagation thesis}) is proved.

Theorem \ref{thm deg} is now completed by Corollary \ref{cor levelset- generalized flow}.
 \end{proof}

\section{The second asymptotic problem}

In this section we consider the asymptotics of the solutions of the singular Cauchy problem for the equation (\ref{allen cahn eps eq})
\begin{equation}\label{asymp cauchy problem allen cahn}
\left\{
\begin{array}{ll}
\mbox{(i)}&u^{\gep}_t-\tr{\left(A^\gep(x)D^2 u^\gep+(\nabla u^\gep(x)D\sigma^\gep(x))\sigma^\gep(x)\right)}+\frac{f(u^\gep)}{\gep^2}=0\quad\mbox{in }\RRn\times(0,+\infty),\\
\mbox{(ii)}&u^{\gep}(x,0)=g(x)\quad\mbox{  in }\RRn,
\end{array}
\right.
\end{equation}
where now the exponent $k$ that appears in (\ref{sigmaeps def}) is a fixed real number so that $0<k<1/4$ and $g$ is a continuous real function in $\RRn$ which takes values in $[m_-,m_+]$ on both sides of $m_o$. Moreover we suppose that $(q(\cdot),0)$ is a solution of the traveling wave equation \eqref{traveling wave equation}. 
Our core assumptions in this section require that $\sigma=(\sigma_{ik})_{i,k}\in C^4(\RRn,\RR^{n\times m})\cap W^{4,\infty}(\RRn,\RR^{n\times m})$ satisfies (\ref{A regularity}), (\ref{A non negative}), $f$ satisfies \eqref{f property}, the traveling wave $q(\cdot)$ satisfies \eqref{q property}, \eqref{q property 2} and the adjoint solution $p(\cdot)$ satisfies \eqref{p convergence}, \eqref{p estimates}.

Our proof uses the definition of generalized flow where the test function satisfies all conditions (i)--(v). 
In this case, we proved in Theorem 2 in \cite{deso5}, combined with Theorem \ref{thm:level set-viscosity sol} that if the comparison principle holds for the subelliptic mean curvature flow equation and $\RRn$ can be endowed with a particular group law so that $(\RRn,\circ)$ is isomorphic to a Carnot group of step two,
then the equivalence of generalized flows and the level set approach is satisfied. This is reported here in Theorem \ref{flow-viscositysol}.
However the proofs in this section will not need this assumption directly, unless when we apply Corollary \ref{cor levelset- generalized flow} at the end of the proof of Theorem \ref{thm2 deg}. 
We will keep the drift term in \eqref{asymp cauchy problem allen cahn}, while it disappears in the step two Carnot group setting, to emphasize that the assumption is not really needed here.
 We define two regions $\Omega^1$ and $\Omega^2$ as in (\ref{omega1-2B}), and show that the front separating $\Omega^1$ and $\Omega^2$ evolves according to the geometric pde
\begin{equation}\label{deg mean curvature eq}
u_t(x,t)+F(x,X u(x,t),X^2u(x,t))=0,\quad(x,t)\in\RRn\times(0,+\infty),
\end{equation}
where 
\begin{equation}\label{deg mean curvature F}
\begin{array}{ll}
F(x,\sigma^t(x)p,\sigma^t(x)B\sigma(x))&=G(x,p,B)\\
&=-\tr\Big[\big(\sigma^T(x)B\sigma(x)+(pD\sigma(x))\sigma(x)\big) \Big(I-\frac{p\sigma(x)\otimes p\sigma(x)}{\abs{p\sigma(x)}^2}\Big)\Big],
\end{array}\end{equation}
and $G:\{(x,p)\in\RRn\times\RRn:p\sigma(x)\neq0\}\times\mathcal{S}^n\to\RR$.

We prove the following result.
\begin{thm}\label{thm2 deg} Suppose that the vector fields $X_1,\dots,X_m$ are generators of a {\it step two Carnot group}.
Assume that the matrix map $A(x)=\sigma(x)\sigma^{t}(x)$ and the function $f$ satisfy (\ref{A regularity}), (\ref{A non negative}) (\ref{f property}). Let $(q(\cdot),0)$ be the solution of the traveling wave equation (\ref{traveling wave equation}), satisfying also (\ref{q property}) and (\ref{q property 2}). 
Finally assume that \eqref{deg mean curvature eq} satisfies the {\it comparison principle}.\\
For $\gep>0$, let $u^\gep$ be the unique solution of the Cauchy problem (\ref{asymp cauchy problem allen cahn}), with $k\in(0,1/4)$ and $g:\RRn\to[m_-,m_+]$ be a continuous function such that the sets $\Gamma_o=\{x:g(x)=m_o\}$, $\Omega_o^+=\{x:g(x)>m_o\}$, $\Omega_o^-=\{x:g(x)<m_o\}$ are nonempty and $(\Gamma_o,\Omega_o^+,\Omega_o^-)\in{\mathcal E}$. Then
$$u^\gep(x,t)\longrightarrow
\left\{\begin{array}{rll}
m_+&\mbox{in} &\{(x,t):v(x,t)>0\},\\
m_-&\mbox{in} &\{(x,t):v(x,t)<0\},
\end{array}
\right.
$$
locally uniformly as $\gep\to0$,
where $v$ is the unique viscosity solution of
\begin{equation}\label{PCfront two}
\left\{\begin{array}{ll}
v_t(x,t)+G(x,\nabla v(x,t),D^2v(x,t))=0\mbox{ in }\RRn\times(0,+\infty),\\
v(x,0)=d_o(x),
\end{array}\right.
\end{equation}
 $G$ is as in (\ref{deg mean curvature F}) and $d_o$ is the signed distance to $\Gamma_o$ which is positive in $\Omega^+_o$ and negative in $\Omega_o^-$. 
\end{thm}

\begin{rem}
We comment here on the relevant differences between the assumptions of Theorem~\ref{PCfront two} and Theorem~\ref{thm deg}, as well as their consequences in the proofs, which may not be immediately obvious:
\begin{itemize}
\item \textbf{Comparison principle:} Assumption that Equation~\eqref{deg mean curvature eq} satisfies the comparison principle. As discussed in the introduction, the general theory for the comparison principle for viscosity solutions of~\eqref{deg mean curvature eq} in the sub-Riemannian framework is not yet complete, but it is required here to ensure a well-posed level-set formulation. Directly, the comparison principle is used only in Step (\emph{~3. Conclusions}) of the general method when invoking Corollary~\ref{cor levelset- generalized flow}, and nowhere else in the calculations of this section.
\item \textbf{Step-two Carnot group structure:} Assumption that the vector fields $X_1, \dots, X_m$ generate a step-two Carnot group. Definition~\ref{generalized flow in Carnot group def} in this section incorporates property~(v) on the test functions, whereas the previous section does not. This has a strong consequence: proving the equivalence between the generalized flow and the level-set approach in the sub-Riemannian framework requires at the moment the assumption that the vector fields generate a step-two Carnot group. This is established not in the present paper, but in our earlier work~\cite{deso5}. Condition~(v) appears similarly in an equivalent definition of viscosity solutions for the mean curvature flow equation, useful for establishing the aforementioned equivalence already in the Euclidean setting, which was introduced by Barles and Georgelin~\cite{bage}, with the equivalence of generalized flows proved by Barles and Souganidis~\cite{bs}. Here, restriction~(v) on the test functions is invoked exclusively when handling term~\eqref{estimate} via Lemma~\ref{lemtech}. It does not appear anywhere else in the computations of this paper, nor does the full geometry of a step-two Carnot group. Except again when invoking Corollary~\ref{cor levelset- generalized flow}. 
\item \textbf{The parameter $k$ in the Riemannian approximation:} In the previous section, it is assumed that $k < 1$, whereas here we require $k < 1/4$. Although we do not seek the optimal choice for $k$, the reader will notice that the core of the construction in the proofs of Lemmas~\ref{deg initialization second lemma2} and~\ref{lemma propagation2} relies on functions~\eqref{eqtest4a} and~\eqref{eqtest4b}, instead of~\eqref{eqtest3a} and~\eqref{eqtest3b}, respectively in Lemmas~\ref{deg initialization second lemma} and~\ref{deg lemma propagation}. Beyond the differences in these test functions (which involve more terms to control), additional derivatives must be computed, yielding this stricter bound. Specifically, the restriction $k < 1/4$ is dictated by a single term $C_4 \varepsilon^{-4k}$ arising from third-order derivatives of $q(\cdot)$ and $\alpha\;p(\cdot)$ near characteristic points (see the bound for~\eqref{III def} and the subsequent formula).
\item \textbf{Higher regularity of the vector fields:} The number of derivatives to be calculated on the functions $q(\cdot)$ and $\alpha\;p(\cdot)$, which are contained in the proposed subsolution functions~\eqref{eqtest4a} and~\eqref{eqtest4b} also explains the greater regularity required of the vector fields in this section, as we compute four derivatives of the vector fields. Notice for instance the term $\alpha(x,t) = \langle Y_\varepsilon\Phi\,Y_\varepsilon^2\Phi,\,Y_\varepsilon\Phi\rangle/|Y_\varepsilon\Phi|^2$, already containing derivatives of the vector fields, and its derivatives appearing in the remainder
estimate~$III_\varepsilon$. 
\end{itemize}
\end{rem}

The proof of the theorem follows along the lines of the previous section, with relevant modifications in the construction of the sub and supersolutions, as discussed in the previous remark. We will concentrate on the main changes needed.

\begin{proof} As in the proof of Theorem \ref{thm deg} we define two families of open sets of $\RRn$ $(\Omega^1_t)_{t\in[0,T)}$ and $(\Omega^2_t)_{t\in[0,T)}$ as in (\ref{omega1-2B}) this time with $\tau=1-k$, (\ref{omega1-2t}), (\ref{omega0}).
By the maximum principle we know that $m_-\leq u^\gep \leq m_+$. 

\emph{First step: initialization.} In this first part we want to show that the traces on $t=0$ satisfy
$$\Omega_o^+=\{d_o>0\}\subseteq\Omega_o^1,\qquad\Omega_o^-=\{d_o<0\}\subseteq\Omega_o^2.$$
We concentrate on the first inclusion.
Let $\hat x$ be such that $d_o(\hat x)>0$, i.e. $g(\hat x)>m_o$. By the continuity of $g$ we can find an $r,\sigma>0$ such that $g(x)\geq m_o+4\sigma,$
for any $x\in B(\hat x,r)$. This means that
\begin{equation}\label{deg lower estimate uepsilonb}
u^\gep(x,0)=g(x)\geq( m_o+4\sigma)\mathds1_{B(\hat x,r)}(x)+m_-\mathds1_{B(\hat x,r)^c}(x).
\end{equation}
Let $\Phi:\RRn\times[0,T]\to\RR$ be defined as in (\ref{deg Phi def}). 
The first Lemma is about the very short time generation of the interface through properties of ordinary differential equations.
We borrow the argument from Theorem 1 of \cite{xce} as revisited in Lemma 2.2 of \cite{bdl}, and in the first parte we only provide the main details for the reader's convenience.
\begin{lem}\label{deg initialization first lemma2} Under the assumptions of Theorem \ref{thm2 deg} we have that for any $\beta>0$ there exist $\tau=\tau(\beta)>0$ and $\bar\gep=\bar\gep(\beta)>0$ such that, for all $0<\gep\leq\bar\gep$, we have
$$u^\gep(x,t_\gep)\geq(m_+-\beta\gep)\mathds1_{\{\Phi(\cdot,0)\geq\beta\}}(x)+m_-\mathds1_{\{\Phi(\cdot,0)<\beta\}}(x),\quad x\in\RRn,$$
where $t_\gep=\tau\gep^2|\lg\gep|$.
\end{lem}
\begin{proof}
Let $\beta>0$ fixed. To prove our thesis we modify the function $f$. 
Consider a smooth cut-off $\rho\in C^{\infty}_0$ such that $0\leq\rho\leq1$ in $\R$, $\rho(s)=1$ if $|s|\leq1$ and $\rho(s)=0$ if $|s|\geq2$. Assume furthermore that $\rho$ satisfies $-2\leq s\rho'(s)\leq0$ and $|\rho''(s)|\leq4$ for all $s\in\R$. 
Now choose $\sigma>0$ so that $[m_o-\sigma,m_o+\sigma]\subset (m_-,m_+)$, define two smooth functions $\rho_1,\rho_2:\R\to[0,1]$ as
$$
\rho_1(q)=\rho\Big(\frac{q-m_o-\sigma/2}{\sigma/2}\Big)\qquad\rho_2(q)=\rho\Big(\frac{q-m_o-\sigma/2}{\sigma/4}\Big)
$$
and set
$$\bar f(q)=(1-\rho_1(q))f(q)+\rho_1(q)f(q-\sigma/2)$$
and
$$\tilde f^\gep(q)=(1-\rho_2(q))\bar f(q)+\rho_2(q)\frac{\frac{\sigma}{2}+m_o-q}{|\lg\gep|}.$$
Notice that $\tilde f^\gep(q)=\overline f(q)=f(q)$ if $q\notin(m_o-\sigma/2,m_o+3\sigma/2)$
and that $\tilde f^\gep(q)=(\sigma/2+m_o-q)/|\log\gep|$ if $q\in [m_o-\sigma/4,m_o+3\sigma/4]$.
Therefore $\tilde f^\gep$ satisfies properties similar to those of $f$ with zeroes $\{m_-,m_o+\sigma/2,m_+\}$. Moreover $f\leq\tilde f^\gep$, for $\sigma$ small.

1. We denote by $\chi=\chi(\tau,\xi)\in C^2([0,+\infty)\times\RRn)$ the solution of
\begin{equation*}
\left\{
\begin{array}{l}
\dot\chi(\tau,\xi)+\tilde f^\gep(\chi(\tau,\xi))=0,\quad \tau>0,\\
\chi(0,\xi)=\xi,
\end{array}
\right.
\end{equation*}
it follows that $\chi_\xi>0$ (as in (\ref{property one chi}) of Lemma \ref{deg initialization second lemma}), while properties (\ref{property two chi}) and (\ref{property three chi}) can be replaced by the following new ones (see \cite{xce}):

\noindent
there exists $\tau_o=\tau_o(\beta,\sigma),\gep_o=\gep_o(\beta,\sigma)>0$ such that, for all $\tau\geq\tau_o|\lg\gep|$ and $\gep\leq\gep_o$
\begin{equation}\tag{$\tilde\chi2$}\label{property two tilde chi}
\begin{array}{l}
\chi(\tau,\xi)\geq m_+-\beta\gep\,\quad\forall\,\xi\geq m_o+\frac{3}{2}\sigma.
\end{array}
\end{equation}
Moreover for any $C\geq\max\{\abs{m_-},\abs{m_+}\}$ we have $\chi(\tau,\xi)\in[-C,C]$ for all $\xi\in[-C,C]$, $\tau\geq0$. Thus,
for any $a>0$ there exists a constant $M_{C,a}>0$ such that
\begin{equation}\tag{$\tilde\chi3$}\label{property three tilde chi}
|\chi_{\xi\xi}(\tau,\xi)|\leq \frac{M_{C,a}}{\gep}\chi_\xi(\tau,\xi),\quad
\end{equation}
for any $\tau\leq a|\ln\gep|$, $\xi\in[-C,C]$, and $\gep$ small enough.

2. Consider a smooth nondecreasing function $\psi:\R\to\R$  such that $\psi(z)=m_-$ if $z\leq0$ and $\psi(z)=m_o+4\sigma$ if $z\geq\beta\wedge\frac{\sigma}{2}$. 
We can prove that the function 
$$\underline u^\gep(x,t)=\chi\Big(\frac{t}{\gep^2},\psi(\Phi(x,0))-\frac{Kt}{\gep}\Big)$$
is a subsolution of (\ref{asymp cauchy problem allen cahn}-i) in $\RRn\times(0,\tau_o\gep^2|\lg\gep|)$ if we choose $K$ large enough.
In fact, from 
$$\begin{array}{cc}\underline u^\gep_t(x,t)=\frac{\dot \chi}{\gep^2}-K\frac{\chi_\xi}{\gep},\quad \nabla\underline u^\gep(x,t)=\chi_\xi(\psi')\nabla\Phi,\\
 D^2\underline u^\gep(x,t)=(\chi_{\xi\xi}(\psi')^2+\chi_\xi(\psi''))\nabla\Phi\otimes\nabla\Phi+\chi_\xi(\psi')D^2\Phi
\end{array}$$
and since $f\leq\tilde f^\gep$, we directly compute, also using \eqref{property three tilde chi}
$$\begin{array}{ll}
\underline u^\gep_t-\tr(Y_\gep\underline u^\gep)+\frac{f(u^\gep)}{\gep^2}&=
\frac{\dot\chi+f(\chi)}{\gep^2}-K\frac{\chi_\xi}{\gep}-(\chi_{\xi\xi}(\psi')^2+\chi_\xi(\psi''))\langle A^\gep\nabla\Phi,\nabla\Phi\rangle\\
&\quad-2\chi_\xi(\psi')\tr(A^\gep)+\chi_\xi(\psi')\tr((\nabla\Phi D\sigma^\gep)\sigma^\gep)
\\ &\leq|\chi_{\xi\xi}|C_1+\frac{\chi_\xi}{\gep}(-K+C_1\gep))\leq 
\frac{\chi_\xi}{\gep}(-K+C_1(M_{C,\tau_0}+\gep))\leq0,
\end{array}
$$
for $K$ large enough.
Furthermore $\underline u^\gep(x,0)\leq u^\gep(x,0)$ by \eqref{deg lower estimate uepsilonb} and the definition of $\psi$. 
Thus using the maximum principle and property (\ref{property two tilde chi}) we can prove that $u^\gep(x,t_\gep)\geq m_+-\beta\gep$ if $t_\gep=\tau_o\gep^2|\lg\gep|$ and $\Phi(x,0)\geq\beta$ for $\gep$ sufficiently small compared to $\beta$. The conclusion then follows.
\end{proof}
Our first main result in this section is the second Lemma.
\begin{lem}\label{deg initialization second lemma2}There exist $\bar h=\bar h(r,\hat x)$, $\bar\beta=\bar\beta(r,\hat x)>0$ independent of $\gep$ such that if $\beta\leq\bar\beta$ and $\gep\leq\bar\gep(\beta)\wedge1$, then there is a subsolution $\omega^{\gep,\beta}$ of (\ref{asymp cauchy problem allen cahn}-i) in $\RRn\times(0,\bar h)$ that satisfies
\begin{equation}\label{eqinitcond}
\omega^{\gep,\beta}(x,0)\leq(m_+-\beta\gep)\mathds1_{\{\Phi(\cdot,0)\geq\beta\}}(x)+m_-\mathds1_{\{\Phi(\cdot,0)<\beta\}}(x),\quad x\in\RRn.\end{equation}
Moreover there is $c>0$ such that if $(x,t)\in B(\hat x,r)\times(0,\bar h)$ and $\Phi(x,t)>3\beta$, then 
$$\liminfstard\frac{\omega^{\gep,\beta}(x,t)-m_+}{\gep^{1-k}}\geq-3c\beta.$$
\end{lem}
\begin{proof}
As in Lemma \ref{deg initialization second lemma} there exist $\gamma,\,\bar h>0$ such that $\bar h<r^2/C$ and
$$\nabla\Phi(x,t)\neq 0,\;Y_\gep\Phi(x,t)\neq0,$$
for any $(x,t)\in Q_{\gamma,\bar h}=\{(x,t):|\Phi(x,t)|\leq\gamma,\,0\leq t\leq\bar h\},\,\gep>0$.
Let $(q(\cdot),0)$ be the solution of the traveling wave equation (\ref{traveling wave equation}) and recall that 
$\Phi(x,t)=r^2-\abs{x-\hat x}^2-Ct$, with $C$ to be chosen sufficiently large.
Suppose that the function $p(r)=\frac{r^2}2\dot q(r)$ is the solution of the inhomogeneous equation (\ref{p equation}) as described in the Appendix, and it satisfies (\ref{p convergence}) and (\ref{p estimates}) and define the smooth function $v^\gep(x,t)$ as
\begin{equation}\label{eqtest4a}
v^{\gep}(x,t)=q\Big(\frac{\Phi(x,t)-2\beta}{\gep\abs{Y_\gep \Phi}}\Big)+\frac{\gep}{\abs{Y_\gep \Phi}}\left(\alpha(x,t)p\Big(\frac{\Phi(x,t)-2\beta}{\gep\abs{Y_\gep \Phi}}\Big)-2\hat c\beta\right),\qquad (x,t)\in Q_{\gamma,\bar h}
\end{equation}
where  
$\alpha(x,t)=\frac{\langle Y_\gep\Phi Y^2_\gep\Phi,Y_\gep\Phi\rangle}{\abs {Y_\gep\Phi}^2}$ and ${\abs{Y_\gep\Phi}}\leq\hat c$ are uniformly bounded with $\gep$ on $Q_{\gamma,\bar h}$. 
We compute
$$v_t^{\gep}(x,t)-\tr(Y_\gep^2 v^\gep(x,t))+\frac{f(v^\gep(x,t))}{\gep^2}=\frac{I_\gep}{\gep^2}+\frac{I\!I_\gep}{\gep}+I\!I\!I_\gep$$
where, similarly as in Lemma \ref{deg initialization first lemma2} we denote $s=\frac{\Phi(x,t)-2\beta}{\gep}$ and also $\lambda=s\abs{Y_\gep \Phi}^{-1}$ to simplify the notation. 
We get
\begin{equation}
\label{I def}I_\gep=-\ddot q+f(q)=0,\end{equation}
\begin{equation}\label{II def}\begin{array}{ll}
I\!I_\gep&=\frac{\dot q}{\abs{Y_\gep \Phi}}\big(\Phi_t-\tr(Y_\gep^2\Phi)\big)-2\Big(\frac{\ddot qs}{\abs{Y_\gep\Phi}}+\dot q\Big)\langle Y_\gep\left(\frac{1}{\abs{Y_\gep\Phi}}\right),Y_\gep\Phi\rangle+(-\ddot p+f'(q)p)\frac{\alpha}{\abs{Y_\gep\Phi}}\\
&\qquad-2\frac{\hat c\beta}{\abs{Y_\gep\Phi}} f'(q)\\
&=\frac{\dot q}{\abs{Y_\gep \Phi}}\big(-C-\tr(Y_\gep^2\Phi)\big)+2\Big(\ddot q(\lambda)\lambda+\dot q(\lambda)\Big)
\frac{\langle Y_\gep\Phi Y_\gep^2\Phi,Y_\gep\Phi\rangle}{\abs {Y_\gep\Phi}^3}+(-\ddot p+f'(q)p)\frac{\alpha}{\abs{Y_\gep\Phi}}\\
&\qquad-2\frac{\hat c\beta}{\abs{Y_\gep\Phi}} f'(q)\\
&=\frac{\dot q}{\abs{Y_\gep \Phi}}\big(-C-\tr(Y_\gep^2\Phi)+\alpha\big)+\frac\alpha{\abs{Y_\gep\Phi}}\Big( 2\ddot q(\lambda)\lambda+\dot q(\lambda) 
-\ddot p+f'(q)p\Big)-2\frac{\hat c\beta}{\abs{Y_\gep\Phi}} f'(q)\\
 &\leq\frac{\dot q(\lambda)}{\abs{Y_\gep\Phi}}\Big(-C-\tr(Y_\gep^2\Phi)+C_3\Big)-2\frac{\hat c\beta}{\abs{Y_\gep\Phi}} f'(q) 
\leq\frac{\dot q(\lambda)}{\abs{Y_\gep\Phi}}\Big(-C+(n+1)C_3\Big)-2\frac{\hat c\beta}{\abs{Y_\gep\Phi}} f'(q)\\
&\leq\frac1{\abs{Y_\gep\Phi}}\left(-\frac C2\dot q(\lambda)-2\hat c\beta f'(q)\right),
\end{array}\end{equation}
by the properties of $p(\cdot)$ in the Appendix and for $C$ sufficiently large.
Furthermore we have the remainder
\begin{multline}\label{III def}
{I\!I\!I}_\gep\leq[2\hat c\gep\beta-\dot qs]\tr\big(Y_\gep^2(\frac{1}{\abs{Y_\gep\Phi}})\big)-\ddot qs^2\big|Y_\gep(\abs{Y_\gep\Phi}^{-1})\big|^2
+\gep\frac\alpha{\abs{Y_\gep\Phi}}[-\dot ps\tr\big(Y_\gep^2(\frac{1}{\abs{Y_\gep\Phi}})\big)\\
-\ddot ps^2\big|Y_\gep(\abs{Y_\gep\Phi}^{-1})\big|^2]
+\alpha\frac{\dot p}{\abs{Y_\gep \Phi}^2}\big(\Phi_t-\tr(Y_\gep^2\Phi)\big)-2\frac\alpha{\abs{Y_\gep\Phi}}\Big(\frac{\ddot ps}{\abs{Y_\gep\Phi}}+\dot p\Big)\langle Y_\gep\left(\frac{1}{\abs{Y_\gep\Phi}}\right),Y_\gep\Phi\rangle\\
+2\gep\langle Y_\gep \frac\alpha{\abs{Y_\gep\Phi}}, Y_\gep p\rangle +\gep p\hbox{Tr}(Y^2_\gep\frac\alpha{\abs{Y_\gep\Phi}})
+\frac{\norm{f_{|[m_-,m_+]}''}_\infty}{2{\abs{Y_\gep\Phi}^2}}(\alpha p-2\hat c\beta)^2.
\end{multline}
Computing higher derivatives, the estimate of $I\!I\!I_\gep$ takes the form
\begin{multline*}
I\!I\!I_\gep\leq [2\hat c\gep\beta-(\dot q+\gep\frac\alpha{\abs{Y_\gep\Phi}}\dot p)s]\tr(Y_\gep^2\big({\abs{Y_\gep\Phi}^{-1}}\big))-(\ddot q+\gep\frac\alpha{\abs{Y_\gep\Phi}}\ddot p)s^2\abs{Y_\gep\big(\abs{Y_\gep\Phi}^{-1}\big)}^2\\
\qquad+\frac{C_3}{\abs{Y_\gep\Phi}^2}|\dot p(\Phi_t-\tr({Y_\gep^2\Phi}))|
+2\left(\frac\alpha{\abs{Y_\gep\Phi}}\right)^2\Big(\frac{\ddot ps}{\abs{Y_\gep\Phi}}+\dot p\Big)+2\dot p\langle Y_\gep \frac\alpha{\abs{Y_\gep\Phi}}, Y_\gep \frac{\Phi-2\hat c\beta}{\abs{Y_\gep\Phi}}\rangle \\
\qquad+\gep p\hbox{Tr}(Y^2_\gep\frac\alpha{\abs{Y_\gep\Phi}})
+\frac{\norm{f_{|[m_-,m_+]}''}_\infty}{2{\abs{Y_\gep\Phi}^2}}(\alpha p-2\hat c\beta)^2
\leq C_4\gep^{-4k},
\end{multline*}
by the notation of Lemma \ref{deg initialization second lemma}, $k<1/4$ and properties (\ref{q property}),  (\ref{p estimates}) of $q$ and $p$, so that $(\abs{\dot q}+\gep^{1-k}\abs{\dot p})\abs{s}$ and $(\abs{\ddot q}+\gep^{1-k}\abs{\ddot p})\abs{s}^2$ are bounded.
Finally
$$
v_t^{\gep}-\tr(Y^2_\gep v^\gep)+\frac{f(v^\gep)}{\gep^2}=\frac{I\!I_\gep}{\gep}+I\!I\!I_\gep
\leq \frac{1}{\gep{\abs{Y_\gep\Phi}}}\Big[-\frac{C}{2}{\dot q(\lambda)}-2\hat c\beta f'(q(\lambda))+C_6\gep^{1-4k}\Big].
$$
To prove that $v^\gep$ is a subsolution of (\ref{asymp cauchy problem allen cahn}-i) in $Q_{\gamma,\bar h}$ it remains to see that the right hand side of the inequality above is non positive. As in Lemma \ref{deg initialization second lemma}, we recall that $f'(m_\pm)\geq d>0$ and $q(r)\to m_\pm$ if $r\to\pm\infty$. Therefore let $M>0$ be a positive constant so that
$$f'(q(r))\geq\frac{d}{2},\quad\mbox{for any }|r|\geq M.$$
Thus, if $\frac{|\Phi-2\beta|}{\gep|Y_\gep\Phi|}\geq M$, we have that
$$v_t^{\gep}(x,t)-\tr(Y_\gep^2 v^\gep(x,t))+\frac{f(v^\gep(x,t))}{\gep^2}\leq\frac{1}{\gep\abs{Y_\gep\Phi}}\Big[-\hat c\;\beta d+C_6\gep^{1-4k}\Big]<-\frac{\beta d}{2\gep}$$
for $\gep$ small enough. If on the contrary $\frac{|\Phi-2\beta|}{\gep|Y_\gep\Phi|}< M$ and we denote with $K$ a strictly positive constant (which depends on $M$) so that  $\dot q(r)\geq K$ for any $|r|\leq M$,  we get that there exists a $\mu=\mu(\beta)$, $0<\mu\leq\beta d/2$ so that, for $\beta$ small compared to $K$, and for small $\gep$,
\begin{equation*}
 v_t^{\gep}-\tr(Y_\gep^2 v^\gep)+\frac{f(v^\gep)}{\gep^2}\leq\frac{1}{\gep\abs{Y_\gep\Phi}}[-\frac{CK}{2}+2\hat c\beta\norm{f_{|[m_-,m_+]}'}+o_\gep(1)]
\leq-\frac{\mu(\beta)}{\gep}<0,
\end{equation*}

2. This gluing part to build a global subsolution, follows closely the second step of Lemma \ref{deg initialization second lemma} in the previous section and Lemma 4.4 in \cite{bs}. Only small adaptations due to the differences in the function $v^\gep$ are needed.
Once we know that $v^\gep$ is a strict subsolution of (\ref{asymp cauchy problem allen cahn}-i) in $Q_{\gamma,\bar h}$ we define, for each $(x,t)\in\{(x,t)\in\RRn\times[0,\bar h]:\Phi(x,t)\leq\gamma\}$,
$$\bar v^\gep(x,t)=
\left\{\begin{array}{lll}
\sup(v^\gep(x,t),m_-)&\mbox{if}&-\gamma<\Phi(x,t)\leq\gamma,\\
m_-&\mbox{if}&\Phi(x,t)\leq-\gamma.
\end{array}
\right.$$
We prove that $\bar v^\gep$ is a continuous viscosity subsolution of equation (\ref{asymp cauchy problem allen cahn}-i) in $\{(x,t)\in\RRn\times[0,\bar h]:\Phi(x,t)\leq\gamma\}$, for $\gep$ sufficiently small. This is obvious in the set $\{\abs{\Phi}\leq\gamma\}$ since $\bar v^\gep$ is the supremum of two subsolutions and if $-\gamma\leq\Phi(x,t)\leq-\gamma/2$, by properties (\ref{q property}) and (\ref{p convergence}), since
$$v^\gep(x,t)\leq  m_-+a e^{-\frac{b(\gamma+4\beta)}{2\gep\hat c}}+\frac\gep{\abs{Y_\gep\Phi}}(C_1p-2\hat c\beta)\leq m_-$$
for small $\gep$ since $p=o_\gep(1)$ in that region, for $\beta$ small.
Thus $\bar v^\gep(x,t)=m_-$. 

Finally we can define $\omega^{\gep,\beta}:\RRn\times[0,\bar h]\to\RR$ as
$$
\omega^{\gep,\beta}(x,t)=
\left\{\begin{array}{ll}
\psi(\Phi(x,t))\bar v^{\gep}(x,t)+(1-\psi(\Phi(x,t)))(m_+-\beta\gep)&\mbox{if }\Phi(x,t)<\gamma,\\
m_+-\beta\gep&\mbox{if }\Phi(x,t)\geq\gamma,
\end{array}
\right.
$$
where $\psi:\RR\to\RR$ is a smooth function such that $\psi'\leq0$ in $\RR$, $\psi=1$ in $(-\infty,\gamma/2]$, $0<\psi<1$ in $(\gamma/2,3\gamma/4)$ and $\psi=0$ in $[3\gamma/4,+\infty)$. The only subset of $\RRn\times(0,\bar h)$ in which we have to check that $\omegaeb$ is a subsolution of (\ref{asymp cauchy problem allen cahn}-i) is $\{(x,t)\in\RRn\times(0,\bar h):\gamma/2\leq \Phi(x,t)\leq3\gamma/4\}$. We have
\begin{multline}\label{computation4}
\omegaeb_t-\tr(Y^2_\gep\omegaeb)+\frac{f(\omegaeb)}{\gep^2}=[\psi'(\Phi_t-\tr(Y^2_\gep\Phi))-\psi''\abs{Y_\gep\Phi}^2] (\bar v^\gep-(m_+-\beta\gep))\\
+\psi(\bar v^\gep_t-\tr(Y^2_\gep\bar v^\gep))-2\psi'\langle Y_\gep\Phi, Y_\gep\bar v^\gep\rangle+\frac{f(\omegaeb)}{\gep^2}
\end{multline}
Take $2\beta<\gamma/4$, then 
$\lambda=(\Phi(x,t)-2\beta)/(\gep\abs{Y_\gep\Phi})>\gamma/(4\hat c\gep)$ and as before
$$v^\gep(x,t)\geq m_+-ae^{\frac{-b\gamma}{4\hat c\gep}}+\frac\gep{\abs{Y_\gep\Phi}}(o_\gep(1)-2\hat c \beta)>m_o,$$
for small $\gep$. Thus $\bar v^\gep(x,t)=v^\gep(x,t)$ and similarly, for $\gep$ small
$$v^\gep(x,t)-(m_+-\beta\gep)\leq \frac\gep{\abs{Y_\gep\Phi}}(o_\gep(1)-\hat c\beta)\leq0.$$
Moreover by (\ref{f property}), $f$ is convex in a neighborhood of $m_+$ and $f(\omegaeb)\leq\psi f(v^\gep)+(1-\psi)f(m_+-\beta\gep)$. The equality (\ref{computation4}) thus becomes, by the end of the first step (here $C$ is large enough),
$$\begin{array}{ll}
\omegaeb_t-\tr(Y^2_\gep\omegaeb)+\frac{f(\omegaeb)}{\gep^2}&\leq[\psi'(-C-\tr(Y^2_\gep\Phi))-\psi''\abs{Y_\gep\Phi}^2] (v^\gep-(m_+-\beta\gep))\\
&\qquad\;-\psi\frac{\mu}{\gep}-2\psi'\langle Y_\gep\Phi, Y_\gep v^\gep\rangle+(1-\psi)\frac{f(m_+-\beta\gep)}{\gep^2}\\
&\leq-\psi\frac{\mu}{\gep}+(1-\psi)\frac{f(m_+-\beta\gep)}{\gep^2}\\
&\qquad-\psi''\abs{Y_\gep\Phi}^2(v^\gep-(m_+-\beta\gep))-2\psi'\langle Y_\gep\Phi, Y_\gep v^\gep\rangle\\
&=-\psi\frac{\mu}{\gep}+(1-\psi)\frac{f(m_+-\beta\gep)}{\gep^2}+C_6(1+\gep^{-k}+\gep^{1-2k}).
\end{array}$$
Finally we observe that for $\gep$ small and by assumption (\ref{f property}) on $f$, (here $0<d/2\leq f'$ near $m_+$)
$$\frac{f(m_+-\beta\gep)}{\gep^2}=\frac{f(m_+-\beta\gep)-f(m_+)}{\gep^2}\leq -\frac{d\beta}{2\gep}<0.$$
Therefore
$$\omegaeb_t-\tr(Y^2_\gep\omegaeb)+\frac{f(\omegaeb)}{\gep^2}\leq\frac{1}{\gep}[-\psi\mu-(1-\psi)\frac{d\beta}2+C_7\gep^{1-k}]<0,$$
for $\gep$ small enough.

3. We now want to examine the request \eqref{eqinitcond} on the initial condition $\omegaeb(\cdot,0)$. By construction $\omega^{\gep,\beta}(x,0)=m_+-\beta\gep$ for $\Phi(x,0)\geq\gamma$, and $\omega^{\gep,\beta}(x,0)=m_-$ for $\Phi(x,0)\leq-\gamma$, hence the required estimate holds automatically outside the compact set $\{x:|\Phi(x,0)|\leq\gamma\}$.
We first observe that, by (\ref{p convergence}), given $\beta$ there exists a $\rho>0$ such that, if $\abs{u}\geq\rho$, then
$$\abs{p(u)}\leq\frac{\beta\hat c}{\|\alpha\|_\infty}.$$
If $(x,t)\in Q_{\gamma,\bar h}$ is such that $\abs{\Phi(x,t)-2\beta}\geq\rho\;\gep\abs{Y_\gep\Phi(x,t)}$, then
\begin{equation}\label{inequality1}
v^\gep(x,t)\leq q\Big(\frac{\Phi(x,t)-2\beta}{\gep\abs{Y_\gep\Phi(x,t)}}\Big)+\frac\gep{\abs{Y_\gep\Phi}}\left(\alpha(x,t)\;p\Big(\frac{\Phi(x,t)-2\beta}{\gep\abs{Y_\gep\Phi(x,t)}}\Big)-2\hat c\beta\right)\leq m_+-\beta\gep.
\end{equation}
If instead $\abs{\Phi(x,t)-2\beta}\leq\rho\;\gep\abs{Y_\gep\Phi(x,t)}$, $(x,t)\in Q_{\gamma,\bar h}$, let $\nu(\rho)>0$ be a constant so that
$$q(u)\leq m_+-\nu(\rho),\qquad \hbox{for any }\abs{u}\leq\rho.$$
Therefore, for $\gep$ small enough,
\begin{equation}\label{inequality2}
v^\gep(x,t)\leq m_+-\nu(\rho)+\frac\gep{\abs{Y_\gep\Phi}}\left(\norm{\alpha}_\infty\norm{p}_\infty-2\hat c\beta\right)\leq m_+-2\beta\gep.
\end{equation}
Combining the estimates in (\ref{inequality1}) and (\ref{inequality2}) for $v^\gep$ we can conclude that $\omegaeb(x,t)\leq m_+-\beta\gep$ for any $(x,t)\in\RRn\times[0,\bar h]$, in particular at $t=0$. 

Assume now that $\Phi(x,t)<\beta$ and $\gep$ is such that $\gep\rho\hat c\leq\beta$.
We obtain $\Phi(x,t)-2\beta<-\beta\leq-\rho\gep\abs{Y_\gep\Phi(x,t)}$, for any $(x,t)\in Q_{\gamma,\bar h}$. Furthermore
$$v^\gep(x,t)\leq q\Big(\frac{-\beta}{\gep\abs{Y_\gep\Phi(x,t)}}\Big)+\frac\gep{\abs{Y_\gep\Phi}}\left(\alpha p\Big(\frac{\Phi-2\beta}{\gep\abs{Y_\gep\Phi(x,t)}}\Big)-2\beta\hat c\right)
\leq m_-+ae^{-\frac{b\beta}{\gep\hat c}}-\beta\gep\leq m_-,$$
if $\gep$ is small.
Therefore $\Phi(x,t)<\beta$ implies $\bar v^\gep(x,t)=m_-$ and, since we have assumed $\beta<\gamma/8$, we conclude that $\omegaeb(x,t)=m_-$. Thus
$$\omegaeb(x,t)\leq(m_+-\beta\gep)\mathds1_{\{\Phi\geq\beta\}}(x,t)+m_-\mathds1_{\{\Phi<\beta\}}(x,t),$$
for any $(x,t)\in\RRn\times[0,\bar h]$, in particular at $t=0$.

To complete the Lemma, with a reasoning similar to the previous one, we can prove that if $(x,t)\in B(\hat x,r)\times(0,\bar h)$, $\gamma\geq\Phi(x,t)>3\beta$, and $\gep$ is small,
  then
$$\omega^{\gep,\beta}(x,t)\geq q\Big(\frac{\beta}{\gep\abs{Y_\gep\Phi}}\Big)-3\hat c\beta\frac\gep{\abs{Y_\gep\Phi}}\geq m_+-ae^{-\frac{b\beta}{\gep\hat c}}-3\beta C_7\gep^{1-k}.$$
Hence $\liminfstard\frac{\omegaeb(x,t)-m_+}{\gep^{1-k}}\geq-3C_7\beta$ for any $(x,t)\in B(\hat x,r)\times(0,\bar h)$ with $\Phi(x,t)>3\beta$.
\end{proof}
In order to end the first step after the two Lemmas, and to prove that $(\Omega_o^+\ni)\hat x\in\Omega^1$, we can now follow the corresponding argument of the one at the end of Lemma \ref{deg initialization second lemma}.

\emph{Second step: propagation.} In this step we show that $\omegaone$ and $((\Omega^2_t)^c)_{t\in(0,T)}$, defined as in \eqref{omega1-2B}, \eqref{omega0}, are respectively super and subflow with normal velocity $-F$ with $F$ defined as in (\ref{deg mean curvature F}). We only show that $\omegaone$ is a superflow.

Let $x_0\in\RRn$, $t\in(0,T)$, $r>0$, $h>0$ so that $t+h<T$, and suppose that $\phi:\RRn\times[0,T]\to\RR$ is a smooth function such that, for a suitable $\tilde C>0$,
$$\phi_t(x,s)+G^*(x,\nabla\phi(x,s),D^2\phi(x,s))\leq-\tilde C<0, \hbox{for all }(x,s)\in B(x_0,r]\times[t,t+h].$$
Assume that $\phi$ also satisfies the properties (ii)--(v) in Definition \ref{generalized flow in Carnot group def}.
We have to show that for every $s\in(t,t+h)$,
$$\{x\in B(x_0,r]:\phi(x,s)>0\}\subset\Omega_s^1.$$
Using the assumptions and the definition of $(\Omega^1_t)_{t\in(0,T)}$ this is equivalent to prove that for all $x\in B(x_0,r)$, $s\in(t,t+h)$ such that $\phi(x,s)>0$, we have
\begin{equation}\label{propagation thesis2}\liminfstard \Big(\frac{u^\gep(y,\tau)-m_+}{\gep^{1-k}}\Big)\geq0\end{equation}
for $(y,\tau)$ in a neighborhood of $(x,s)$.  
As in Theorem \ref{thm deg}, the main part of the proof is the following lemma.

\begin{lem}\label{lemma propagation2} Let $x_o,t,r,h$ be as above, and let $\phi$ be a smooth function satisfying (i)--(v) in Definition \ref{generalized flow in Carnot group def}.
There exists $\bar\beta$ small enough such that, if $\beta\leq\bar\beta$ and $\gep\leq\bar\gep(\beta)$ then there is a viscosity subsolution $\omegaeb$ of (\ref{asymp cauchy problem allen cahn}-i) in $B(x_0,r)\times(t,t+h)$ that satisfies,

 1. $\omegaeb(x,t)\leq(m_+-\beta\gep)\mathds1_{\{\phi(\cdot,t)\geq\beta\}}(x)+m_-\mathds1_{\{\phi(\cdot,t)<\beta\}}(x),\quad\mbox{for all } x\in B(x_0,r],$

2. $\omega^{\gep,\beta}(x,s)\leq(m_+-\beta\gep)\mathds1_{\{\phi(\cdot,s)\geq\beta\}}(x)+m_-\mathds1_{\{\phi(\cdot,s)<\beta\}}(x),\quad\mbox{for all } x\in \partial B(x_0,r]$, $s\in[t,t+h]$

3. there is $ c>0$ such that if $(x,s)\in B(x_0,r]\times[t,t+h]$ satisfies $\phi(x,s)>3\beta$, then
$$\liminfstard\frac{\omegaeb(x,s)-m_+}{\gep^{1-k}}\geq-3c\;\beta.$$
\end{lem}

\begin{proof}[Proof of Lemma \ref{lemma propagation2}]
This proof is similar to the one of Lemma \ref{deg initialization second lemma2} and we basically exchange $\Phi$ with $\phi$ which has a more general time dependence. We insist on the main differences.
Observe that since $\phi$ satisfies property (ii) in the definition, there exists $\gamma>0$ such that $|\nabla\phi(x,s)|\neq0$ and
$|Y_\gep\phi(x,s)|\neq0$ for any $(x,s)\in Q_\gamma=\{(x,s)\in B(x_0,r]\times[t,t+h]:|\phi(x,s)|\leq\gamma\}$ and all $\gep>0$. 
Let $(q(\cdot),0)$ be the solution of the traveling wave equation and $p(\cdot),\alpha(x,s),\hat c$ be as in the proof of Lemma \ref{deg initialization second lemma2}. We define as our candidate subsolution the function
\begin{equation}\label{eqtest4b}
v^{\gep}(x,s)=q\Big(\frac{\phi(x,s)-2\beta}{\gep|Y_\gep\phi|}\Big)+\frac\gep{\abs{Y_\gep\phi}}\left(\alpha(x,s)p\Big(\frac{\phi(x,s)-2\beta}{\gep\abs{Y_\gep\phi}}\Big)-2\hat c\beta\right),\qquad(x,s)\in Q_\gamma.\end{equation}
It turns out that
$$v^\gep_t(x,s)-\tr(Y^2_\gep v^\gep(x,s))+\frac{f(v^\gep(x,s))}{\gep^2}=\frac{I_\gep}{\gep^2}+\frac{I\!I_\gep}{\gep}+{I\!I\!I}_\gep,$$
where $I_\gep$, $I\!I_\gep$, ${I\!I\!I}_\gep$ are exactly as in (\ref{I def}), (\ref{II def}) and (\ref{III def}), with the function $\Phi$ replaced by $\phi$. 
We get $I_\gep= \ddot q(\lambda)+f(q(\lambda))=0$
and
\begin{equation}\label{computation5}\begin{array}{ll}
I\!I_\gep&=\frac{\dot q(\lambda)}{|Y_\gep\phi|}\Big(\phi_t-\tr(Y^2_\gep\phi)+ \frac{\langle(Y_\gep^2\phi)Y_\gep\phi,Y_\gep\phi\rangle}{\abs{Y_\gep\phi}^2}\Big)-2\frac{\hat c}{\abs{Y_\gep\phi}}\beta f'(q(\lambda))\\
 &\leq\!\frac{\dot q(\lambda)}{|Y_\gep\phi|}\Big(-\tilde C\!-\!G^*(\nabla\phi,D^2\phi)\!-\!\tr(Y_\gep^2\phi\!)\!+\! \frac{\langle Y_\gep\phi(Y_\gep^2\phi),Y_\gep\phi\rangle}{\abs{Y_\gep\phi}^2}\Big)-2\frac{\hat c\beta}{\abs{Y_\gep\phi}} f'(q(\lambda))\\
 &\leq\!\frac{\dot q(\lambda)}{|Y_\gep\phi|}\Big(\!\!-\!\tilde C-\!\gep^{2k}\tr(D^2\phi)\!+\!\frac{\langle Y_\gep\phi Y_\gep^2\phi,Y_\gep\phi\rangle}{\abs{Y_\gep\phi}^2} -\!\langle \frac{X\phi}{{\abs{X\phi}}}X^2\phi ,\frac{X\phi}{\abs{X\phi}}\rangle^*\Big)-2\frac{\hat c\beta}{\abs{Y_\gep\phi}} f'(q(\lambda)),
 \end{array}\end{equation}
where $\lambda=\frac{\phi(x,s)-2\beta}{\gep|Y_\gep\phi(x,s)|}$.
We immediately observe that if $\abs{X\phi(x,s)}=0$, by the assumption (v) on $\phi$ we have $\abs{X^2\phi(x,s)}=0$ and $Y_\gep^2\phi\to0$ uniformly on the subset of such points. Thus
$${I\!I}_\gep\leq\frac{\dot q(\lambda)}{|Y_\gep\phi|}\big(-\tilde C+o_\gep(1)\big)-\frac{2\hat c}{\abs{Y_\gep\phi}}\beta f'(q(\lambda))
\leq\frac1{|Y_\gep\phi|}\left(-\frac{\tilde C\dot q(\lambda)}{2}-2{\hat c\beta} f'(q(\lambda))\right)$$
for $\gep$ small enough.

We want a similar estimate if instead $\abs{X\phi(x,s)}\neq0$. Then since $|Y_\gep\phi|\geq\max\{|X\phi|,|\nabla\phi| \gep^k\}>0$, we obtain
 $$\begin{array}{ll}
 \frac{\langle Y_\gep\phi(Y_\gep^2\phi) ,Y_\gep\phi\rangle}{\abs{Y_\gep\phi}^2}\!\!\!&= \frac{\langle X\phi(X^2\phi),X\phi\rangle}{\abs{Y_\gep\phi}^2}
 +\gep^{2k}\left(\frac{\langle \nabla\phi(D^2\phi)\sigma(x), X\phi\rangle}{\abs{Y_\gep\phi}^2}
 +\frac{\langle X\phi D(\nabla\phi(x)\sigma(x)) ,\nabla \phi\rangle}{\abs{Y_\gep\phi}^2}
 +\gep^{2k}\frac{\langle \nabla\phi\;D^2\phi ,\nabla\phi\rangle}{\abs{Y_\gep\phi}^2}\right)\\
&\leq \frac{\langle X\phi(X^2\phi),X\phi\rangle}{\abs{X\phi}^2} (1-\frac{\gep^{2k}\abs{\nabla\phi}^2}{\abs{Y_\gep\phi}^2})
+C_9\gep^k
 \end{array}$$
 the inequality in (\ref{computation5}) becomes
 \begin{equation}\label{computation8}
 {I\!I}_\gep
\leq\frac{\dot q(\lambda)}{|Y_\gep\phi|}\Big(-\tilde C+\frac{\gep^{2k}\abs{\nabla\phi}^2}{\abs{Y_\gep\phi}^2}\norm{X^2\phi} +C_9\gep^k \Big)-2\frac{\hat c\beta}{\abs{Y_\gep\phi}} f'(q(\lambda)).
\end{equation}
We want to prove that
\begin{equation}\label{estimate}
\frac{\gep^{2k}\abs{\nabla\phi(x,s)}^2}{\abs{X\phi(x,s)}^2+\gep^{2k}\abs{\nabla\phi(x,s)}^2}\norm{X^2\phi(x,s)}=o_\gep(1),
\end{equation}
uniformly for $(x,s)\in Q_\gamma$.
We separate two cases. If for given $\gep>0$ and for $(x,s)\in Q_\gamma$ we have $\abs{X\phi(x,s)}^2>\gep^k$, then
$$\frac{\gep^{2k}\abs{\nabla\phi(x,s)}^2}{\abs{Y_\gep\phi(x,s)}^2}\norm{X^2\phi(x,s)}<\gep^k
\abs{\nabla\phi(x,s)}^2\norm{X^2\phi(x,s)}\leq C_{10}\gep^k.$$
If instead
$\abs{X\phi(x,s)}^2\leq\gep^k$, observe that
$$
\norm{X^2\phi(x,s)}\frac{\gep^{2k}\abs{\nabla\phi(x,s)}^2}{\abs{Y_\gep\phi(x,s)}^2} \leq
\sup\{\norm{X^2\phi(y,r)}:(y,r)\in Q_\gamma,\;\abs{X\phi(y,r)}^2\leq\gep^k\}=o_\gep(1),
$$
as we conclude by using property (v) of $\phi$, and the following elementary Lemma with $f(y,s)=|X\phi(x,s)|^2$ and $g(y,s)=\norm{X^2\phi(y,s)}$.
\begin{lem}\label{lemtech} Let $K\subset\RRn$ be a compact set and $f,g:K\to[0,+\infty)$ be two continuous functions such that
$$\mbox{if $f(x)=0$ for some $x\in K$, then $g(x)=0$.}$$
Then, it holds
$$\lim_{\gep\to0^+}\sup_{x\in K}\{g(x):f(x)\leq\gep\}=0$$
\end{lem}
\begin{proof} Assume by contradiction that there is $c>0$ such that for any $\gep>0$ we can find $x_\gep=x_\gep(c)\in K$ so that
\begin{equation}\label{contradiction}
g(x_\gep)\geq c\qquad\mbox{ and }\qquad f(x_\gep)\leq\gep.
\end{equation}
 By compactness, there is a subsequence $x_{\gep_n}\to\hat x\in K$. Passing to the limit as $\gep_n\to0^+$ in (\ref{contradiction}) we obtain the contradiction $g(\hat x)\geq c>0$ and $f(\hat x)=0$.
\end{proof}
So far we thus proved that
$${I\!I}_\gep\leq\frac{1}{\abs{Y_\gep\phi}}\left(-\frac{\tilde C\dot q(\lambda)}{2}-2\hat c\;\beta f'(q(\lambda))\right)$$
for any $(x,s)\in Q_\gamma$ and for $\gep$ small enough. 

We are left with the term $I\!I\!I_\gep$. 
Here we follow along the proof of the first part of Lemma \ref{deg initialization first lemma}. 
When computing $\partial_t v^\gep$, with respect to \eqref{III def} we have the extra terms
\begin{multline}
\left(\frac\gep{\abs{Y_\gep\phi}}(p\;\partial_t\alpha)+\frac\alpha{\abs{Y_\gep\phi}^2}\dot p\;\partial_t\phi\right)+
\left((\dot q(\lambda)+\frac\gep{\abs{Y_\gep\phi}}\alpha\dot p(\lambda))\frac{\phi-2\beta}{\gep}+\gep(\alpha\;p-2\hat c\beta)\right)\partial_t(\abs{Y_\gep\phi}^{-1})\\
\leq C_{10}(\gep^{1-3k}+\gep^{-2k})+C_{10}\big(\abs{\dot q(\lambda)\lambda}+\gep^{1-k}\abs{\dot p(\lambda)\lambda}+\gep^{1-k}\big)
\frac{\langle Y_\gep\phi, Y_\gep\phi_t \rangle}{\abs{Y_\gep\phi}^2}\\
\leq C_{11}(\gep^{1-3k}+\gep^{-2k}+\gep^{-k}+\gep^{1-2k})\; ,
\end{multline}
by using properties (\ref{q property}) and (\ref{p estimates}) again, and as before
$s=\frac{\phi(x,s)-2\beta}{\gep}$, $\lambda=\frac s{\abs{Y_\gep\phi(x,s)}}$.\\
We conclude as in Lemmas \ref{deg initialization second lemma2} and \ref{deg initialization second lemma}, since $0<k<1/4$, and
$$v_t^{\gep}(x,s)-\tr(Y^2_\gep v^\gep(x,s))+\frac{f(v^\gep(x,s))}{\gep^2}\leq\frac{1}{\gep\abs{Y_\gep\phi}}\Big[-\frac{\tilde C\dot q(\lambda)}{2}+C_{12}\gep^{1-2k}-2\hat c\beta f'(q(\lambda))\Big],$$
for all $(x,s)\in Q_{\gamma}$.

The proof of the extension of $v^\gep$ to a global subsolution of (\ref{asymp cauchy problem allen cahn}-i) $\omega^{\gep,\beta}$ in $B(x_o,r]\times[t,t+h]$ and the proof that such a function satisfies 1, 2, 3 in the statement, follows again closely \cite{bs} and it is similar to that of Lemma \ref{deg initialization second lemma2} (see also Lemma \eqref{deg initialization second lemma}) and we omit it.
\end{proof}
Once we have proved the first two steps (initialization and propagation of the front) the proof of Theorem \ref{thm2 deg} is completed using Corollary \ref{cor levelset- generalized flow} as in Section 3. This last step however contains the reference to the level set approach and its comparison principle. Therefore it is important to justify the assumption. The comparison principle for the mean curvature flow equation was not used elsewhere in this section.
\end{proof}

\appendix
\section{Appendix}

This appendix collects auxiliary results used in the paper. We stress that the material in
Appendix~A.1 (properties of the traveling wave equation and of the auxiliary function
$p(\cdot)$) is used directly in the asymptotic constructions of Sections~3 and~4, and is therefore part of the core argument. 
By contrast, the material in Appendix~A.2 (step-two Carnot groups) and Appendix~A.3 (the level-set
approach and the equivalence with generalized flows) is used only in Step~3 (Conclusion) of
Section~4, i.e.\ only for the identification of the generalized flow with the level-set evolution
in the case $c=0$, and is not needed for the unconditional part of the construction.

\subsection{The traveling wave equation}
A classical tool to construct solutions of the reaction-diffusion equation is the traveling wave equation.
By the properties of $f$ in (\ref{f property}) it can be shown (see for example \cite{aw} and \cite{fml}) that there exists a unique pair $(q(\cdot),c)$, solution of the traveling wave boundary value problem
\begin{equation}
\ddot q+c\dot q=f(q),\qquad \lim_{r\to\pm\infty}q(r)=m_{\pm},\qquad q(0)=m_o.
\end{equation}
We assume that $q$ satisfies the following properties
\begin{equation}\label{q property}
\begin{array}{l}
\dot q\geq0,\qquad\sup_{r\in\RR}[(1+|r|+|r|^2+|r|^3)\dot q(r)+(|r|+|r|^2+|r|^3)\abs{\ddot q(r)}]<+\infty\\
\exists\: a,b>0\mbox{ such that }\\
\begin{array}{l}
m_+-(a/2)e^{-2br}\geq q(r)\geq m_+-ae^{-br},\mbox{ as }r\to+\infty,\\
m_-+(a/2)e^{-2br}\leq q(r)\leq m_-+ae^{-br},\mbox{ as }r\to-\infty,\end{array}
\end{array}
\end{equation}
and
\begin{equation}\label{q property 2}
\mbox{there exists a constant }N>0\mbox{ such that }|\ddot q|\leq N\dot q.
\end{equation}
A typical example of function $f$ where all assumptions can be directly checked is
when $f$ is as in (\ref{f example}). We can compute explicitly ($c$ will only depend on $f$)
$$c=2m_o-m_+-m_-,\qquad q(r)=m_-+\frac{m_+-m_-}{1+e^{-(m_+-m_-)(r+\bar r)}},$$
where $\bar r=\frac{1}{m_+-m_-}\ln\big(\frac{m_o-m_-}{m_+-m_o}\big)$. The constant $N$ in (\ref{q property 2}) turns out to be smaller than $m_+-m_-$.

When $c=0$, following \cite{bs}, we also consider the linearized around $q$, inhomogeneous adjoint
\begin{equation}\label{p equation}
\mathcal -\ddot p+ f'(q)p=\chi
\end{equation}
for functions $\chi\in C(\R^n)\cap L^2(\RR)$ such that 
\begin{equation}\label{integral condition chi}\int_{-\infty}^{+\infty}\chi(s)\dot q(s)ds=0.\end{equation}
Indeed there exists a unique $p\in C^2(\RR)\cap H^1(\RR)$, solution of the equation (\ref{p equation}) 
and it is given by
$$p(s)=\dot q(s)\int^s_{0}\frac{1}{(\dot q(\tau))^2}\Big(\int^\tau_{-\infty}\chi(\eta)\dot q(\eta)d\eta\Big)d\tau.$$
It also satisfies
\begin{equation}\label{p convergence}
sp(s)\to0,\qquad\mbox{as }\abs{s}\to+\infty
\end{equation}
and
\begin{equation}\label{p estimates}
\sup_{s\in\RR}(\abs{\dot p(s)}+(1+\abs{s})\abs{\ddot p(s)})<+\infty.
\end{equation}
Notice that in the particular case that $\chi(r)=-(2r\ddot q(r)+\dot q(r))$, then $p(r)=\frac{r^2}2 \dot q(r)$ as directly checked.

In the special case of $f$ as in \eqref{f example} with $m_-=-1$, $m_o=0$, $m_+=1$, the nonlinearity is $f(q)=2q(q^2-1)$. The traveling wave equation (\ref{traveling wave equation}) then simplifies to
$$\ddot q=f(q)$$
and its explicit solution is given by $q(s)=\tanh(s)$. Hence all the properties above for $p(\cdot)$ are satisfied.

\subsection{Step two Carnot groups}

Following Section 4, in order to complete in \cite{deso5} the equivalence of generalized flows with the level set approach, the $m$ vector fields in $\mathcal{X}$ are assumed to be 
generators of a step-two Carnot group. Accordingly, the following H\"{o}rmander condition is satisfied
$$
\operatorname{span}\{X_i(x),\, [X_j, X_k](x) : i, j, k \in \{1, \dots, m\}\} 
= \mathbb{R}^n, \quad \text{for all } x \in \mathbb{R}^n.
$$
Moreover, the space $\mathbb{R}^n$ of a Carnot group is endowed with a group 
structure via the operation
$$
x \circ y = \bigl(x_h + y_h,\; x_v + y_v + \langle Bx_h, y_h \rangle\bigr),
$$
where $B^{(j)}$, $j \in \{1, \dots, n-m\}$, are skew-symmetric, linearly independent $m \times m$ matrices, $Bx_h = (B^{(1)}x_h, \dots, B^{(n-m)}x_h)$, and 
$\langle Bx_h, y_h \rangle = (B^{(1)}x_h \cdot y_h, \dots, 
B^{(n-m)}x_h \cdot y_h)$.
Under this operation, the identity element is $0$ and the inverse is $x^{-1} = -x$.
In this group setting, it is not restrictive to assume that we can decompose 
$x = (x_h, x_v) \in \mathbb{R}^n = \mathbb{R}^m \times \mathbb{R}^{n-m}$, 
and
\begin{equation}\label{smatrix}
\sigma(x) = \begin{pmatrix} I_m \\ (Bx_h)^T \end{pmatrix},
\end{equation}
where $I_m$ is the $m \times m$ identity matrix. With this choice $\sigma$ only depends on the first (horizontal) $m$ coordinates and
given a smooth function $u \in C^2(\mathbb{R}^n)$, we then get as before the 
\emph{horizontal gradient} (expressed as a row vector) by
$Xu(x) = \nabla u(x)\, \sigma(x),$
and the \emph{horizontal Hessian} by
$$
X^2 u(x) = \bigl(X_j X_k\, u(x)\bigr)^*_{j,k=1,\dots,m} 
= \sigma^T(x)\, D^2 u(x)\, \sigma(x).
$$
Notice that the drift term disappeared. Prototypical examples of step two Carnot groups are the Heisenberg groups in $\R^{2n+1}$.

All of the above is classical, and we refer the reader to the monograph by 
Bonfiglioli, Lanconelli, and Uguzzoni \cite{bolaug} for a thorough introduction to the subject.

\subsection{The level set approach and generalized evolution}
In this paragraph, we outline the relation of generalized flows, as defined in Section 2 and the level set approach for geometric equations.
For the details on the level set approach, the reader may consult Souganidis \cite{soug} or Giga \cite{gi}.
Let $\mathcal{E}$ be the collection of all the triples $(\Gamma_0, \Omega^+_0, \Omega^-_0)$ of mutually disjoint subsets of $\RRn$ such that $\Gamma_o$ is closed, $\Omega^{\pm}_o$ are open and $\RRn=\Gamma_0\cup \Omega^+_0\cup \Omega^-_0$. 
We recall the level set approach to front evolution, as linked to the following Cauchy problem
\begin{equation}\label{level set equation}\left\{\begin{array}{ll}
v_t(x,t)+G(x,\nabla v(x,t),D^2v(x,t)),&(x,t)\in\RRn\times(0,T),\\
v(x,0)=v_o(x),&x\in\RRn,
\end{array}
\right.
\end{equation}
where $G$ is a geometric hamiltonian.
\begin{defin}\label{deflevelset} Let $(\Gamma_o,\Omega^+_o,\Omega^-_o)\in{\mathcal E}$ and $\Gamma_o\neq\emptyset$. 
Assume that $G$ satisfies properties (F1-2-3) of Section 2.
A collection of triplets $(\Gamma_t,\Omega_t^+,\Omega^-_t)_{t\in[0,T)}\subset {\mathcal E}$
is a \emph{ level set evolution} of $(\Gamma_o,\Omega_o^+,\Omega^-_o)$ with normal velocity $V=-G$, 
if there exists a viscosity solution of (\ref{level set equation}) such that, for any $t\geq0$,
$$\Gamma_t=\{x\in\RRn:v(x,t)=0\},\quad\Omega_t^+=\{x\in\RRn:v(x,t)>0\}\quad\Omega_t^-=\{x\in\RRn:v(x,t)<0\},$$
\end{defin}
It is well known that the {\it comparison principle} implies that the definition is well posed for sets, i.e. the zero-level set of the solution only depends on the zero-level set of the initial condition.
One of the topological peculiarities of the level set approach is the possibility of the so called {\it fattening} of the fronts. We avoid it if the following definition holds.
\begin{defin}\label{no interior condition} Let $(\Gamma_t,\Omega_t^+,\Omega_t^-)_{t\in[0,T)}\subset{\mathcal E}$ be the level set evolution.
We say that the \emph{no-interior condition } holds for $(\Gamma_t)_{t\in[0,T)}$ if
\begin{equation}
\begin{array}{c}
\closure\{(x,t):v(x,t)>0\}=\{(x,t):v(x,t)\geq0\}\;\mbox{and}\\
\interior\{(x,t):v(x,t)\geq0\}=\{(x,t):v(x,t)>0\},
\end{array}
\end{equation}
where $v$ is a solution of (\ref{level set equation}) as in Definition \ref{deflevelset}.
\end{defin}
Clearly if the no-interior condition in Definition \ref{no interior condition} holds then the set $\cup_{t\geq0}(\Gamma_t\times\{t\})$ has an empty interior in $\RRn\times[0,+\infty)$. 
The following also holds. It explains consequences of the fattening problem and introduces discontinuous solutions.
\begin{thm}\label{thm:level set-viscosity sol} Suppose that $G$ satisfies (F1-2-3), and let $(\Gamma_t,\Omega^+_t,\Omega^-_t)\in{\mathcal E}$, $t\in[0,T]$, be a level set evolution.
\begin{description}
\item[(i)] The two functions $\overline\chi(x,t)=\mathds1_{\Omega_t^+\cup\Gamma_t}(x)-\mathds1_{\Omega_t^-}(x),$
$\underline\chi(x,t)=\mathds1_{\Omega_t^+}(x)-\mathds1_{\Omega_t^-\cup\Gamma_t}(x)$
are (discontinuous) viscosity solutions of (\ref{level set equation}) respectively with the discontinuous initial data $\overline\chi(x,0)=\mathds1_{\Omega_o^+\cup\Gamma_o}-\mathds1_{\Omega_o^-}$ and $\underline\chi(x,0)=\mathds1_{\Omega_o^+}-\mathds1_{\Omega_o^-\cup\Gamma_o}$.
They are the maximal subsolution and the minimal supersolution, respectively.
 \item[(ii)] Suppose that $\Gamma_o$ has an empty interior and that a {\it comparison principle} holds for the equation (\ref{level set equation}). Then the Cauchy problem for (\ref{level set equation}) associated with the initial data $w_o=\mathds1_{\Omega_o^+}-\mathds1_{\Omega_o^-}$ has a unique discontinuous solution if and only if the no-interior condition of Definition \ref{no interior condition} holds. Such solution is given by the function
 \begin{equation}\label{eqchi}
 \chi(x,t)=\mathds1_{\Omega_t^+}(x)-\mathds1_{\Omega_t^-}(x).\end{equation}
\end{description}
\end{thm}
In the previous statement, a locally bounded function $v$ is a discontinuous solution of \eqref{level set equation} if $v^*$ is a supersolution and $v_*$ is a subsolution. Uniqueness of two discontinuous solutions means that they have the same semicontinuous envelopes.
Discontinuous solutions also play a role when we define the \emph{generalized flow} as in Section 2 through the definition of sub-/super-flows, introduced in \cite{bs}, see also \cite{soug}, and then reformulated (and partially revisited) in \cite{bdl}. It can be proved to be equivalent to the level set approach under appropriate assumptions.
\begin{thm} \label{flow-viscositysol}
Assume that the Hamiltonian $G$ in \eqref{level set equation} has the first order structure \eqref{asymp eq} or the second order structure \eqref{allen cahn eps eq}.
Define the generalized flow by Definition \ref{generalized flow in Carnot group def} with properties (i)--(iv) of the test functions in the first case, properties (i)--(v) in the second. Moreover in the latter also assume that the family of vector fields gives to $\R^n$ the structure of a step two Carnot group.
The following equivalences between generalized flows and viscosity solutions of \eqref{level set equation} hold true.
\begin{description}
\item[(i)]Let $\open$ be a family of open subsets of $\RRn$ such that the set $\Omega$, $\Omega:=\bigcup_{t\in(0,T)}\Omega_t\times\{t\}$, is open in $\RRn\times[0,T]$. Then $\open$ is a generalized superflow with normal velocity $V=-G$ in \eqref{level set equation} if and only if the function ${\underline w}=\mathds{1}_{\Omega}-\mathds{1}_{\Omega^c}$ is a viscosity supersolution of (\ref{level set equation}).
\item[(ii)]Let $\close$ be a family of close subsets of $\RRn$ such that the set $\mathcal{F}:=\bigcup_{t\in(0,T)}\mathcal{F}_t\times\{t\}$ is closed in $\RRn\times[0,T]$. Then $\close$ is a generalized subflow with normal velocity $V$ if and only if the function $\overline w=\mathds{1}_{\mathcal{F}}-\mathds{1}_{\mathcal{F}^c}$ is a viscosity subsolution of (\ref{level set equation}).
\end{description}
\end{thm}
\begin{rem}
The previous result is proved in \cite{bdl} in the classical cases and uses only properties (i)--(iv) of the definition of generalized flow. It extends to our first order case basically unchanged, i.e. when $G(x,p)=c|p\sigma(x)|$. In the second order case, when $G$ has the structure of the subelliptic mean curvature operator, we proved it in \cite{deso5} under the stricter conditions (i)--(v) on the test functions, when the vector fields generate a step two Carnot group.
\end{rem}
Putting things together we finally obtain the following result which is crucial in the conclusion step of the abstract method for studying singular limits of reaction diffusion equations. The general idea is from \cite{bs} later revisited in \cite{bdl}.
\begin{cor}\label{cor levelset- generalized flow}
Suppose that the geometric equation in \eqref{level set equation} satisfies the {\it comparison principle} between super and subsolutions.
 Assume we have two families of open subsets of $\RRn$, $(\Omega^1_t)_{t\in(0,T)}$ and $(\Omega^2_t)_{t\in(0,T)}$ such that $(\Omega^1_t)_{t\in(0,T)}$ and $((\Omega^2_t)^c)_{t\in(0,T)}$ are respectively super- and subflows with normal velocity $-G$ and also $\cup_{t\in(0,T)}\Omega^1_t\times\{t\}$, $\cup_{t\in(0,T)}\Omega^2_t\times\{t\}$ are open and disjoint. 
 Define $\underline w=\mathds1_{\Omega^1}(x)-\mathds1_{(\Omega^1)^c}(x)$ and $\overline w=\mathds1_{(\Omega^2)^c}(x)-\mathds1_{\Omega^2}(x)$ and extend them at $t=0$ by lower semicontinuity or upper semicontinuity, respectively.
 
 \noindent
 Suppose moreover that there exists $(\partial\Omega_0^+,\Omega_0^+,\Omega_0^-)\in\mathcal{E}$ such that $\Omega^+_0\subseteq\Omega_0^1$ and $\Omega^-_0\subseteq\Omega_0^2$. Then, if we denote with $(\Gamma_t,\Omega_t^+,\Omega_t^-)_{t\in(0,T)}$ the level set evolution of $(\partial\Omega_0^+,\Omega_0^+,\Omega_0^-)$, we have:\\
\indent (i) for all $t\in[0,T)$,
 $$\Omega^+_t\subset\Omega^1_t\subset\Omega^+_t\cup\Gamma_t,\qquad\Omega^-_t\subset\Omega^2_t\subset\Omega^-_t\cup\Gamma_t,$$ \indent (ii) if $\bigcup_t \Gamma_t\times\{t\}$ satisfies the no-interior condition (\ref{no interior condition}), then for all $t\in[0,T)$, $$\Omega^+_t=\Omega^1_t,\qquad\Omega^-_t=\Omega^2_t.$$
\end{cor}
\begin{proof}
Define $\underline \chi$ and $\overline \chi$ as in Theorem \ref{thm:level set-viscosity sol}. For the first part of the statement of Theorem \ref{flow-viscositysol} (i), $\underline w$ is a supersolution of \eqref{level set equation}. Since $\underline w(x,0)\geq \underline \chi(x,0)$ and by Theorem \ref{thm:level set-viscosity sol} (i), we get that $\underline w\geq \underline \chi$ which is the minimal subsolution. Similarly $\overline w\leq\overline \chi$. The first statement then follows. The second statement is a direct consequence of Theorem \ref{thm:level set-viscosity sol}(ii).
\end{proof}


\end{document}